\documentclass[11pt]{article}

\usepackage[T1]{fontenc}
\usepackage[utf8]{inputenc}
\usepackage[english]{babel}
\usepackage[a4paper,margin=2.75cm]{geometry}

\usepackage{amsmath,amsthm,amssymb,amsfonts}
\usepackage{mathtools}
\usepackage{enumitem}
\usepackage{mathpazo}
\usepackage[scaled=0.92]{helvet}
\usepackage{mathrsfs}
\usepackage{microtype}
\usepackage[cal=boondoxo,scr=boondox]{mathalfa}
\usepackage{xcolor}
\definecolor{paperblue}{HTML}{1E4B66}
\usepackage{titlesec}
\usepackage{titling}
\usepackage{abstract}

\usepackage{nameref}
\usepackage{hyperref}
\usepackage{cleveref}
\hypersetup{
	colorlinks=true,
	linkcolor=paperblue,
	citecolor=paperblue,
	urlcolor=paperblue,
	pdfencoding=auto,
	pdftitle={Codimension-Two Spacelike Submanifolds with Umbilical Lightlike Normal Sections and Their Relationship to Lightlike Hypersurfaces},
	pdfauthor={Juan S. Gómez},
	pdfsubject={Lorentzian geometry and lightlike hypersurfaces},
	pdfkeywords={spacelike submanifolds, lightlike normals, lightlike hypersurfaces, umbilicity, Lorentzian geometry}
}
\pretitle{\begin{center}\LARGE\sffamily\bfseries\color{paperblue}}
\posttitle{\par\end{center}\vspace{0.6em}}
\preauthor{%
	\begin{center}\normalsize\lineskip .45em
	\begin{tabular}[t]{c}}
\postauthor{\end{tabular}\par\end{center}\vspace{-0.8em}}
\predate{}
\postdate{}

\titleformat{\section}
	{\normalfont\Large\sffamily\bfseries\color{paperblue}}
	{\thesection}{0.7em}{}
\titleformat{\subsection}
	{\normalfont\large\sffamily\bfseries\color{paperblue}}
	{\thesubsection}{0.7em}{}
\titlespacing*{\section}
	{0pt}{2.8ex plus .8ex minus .2ex}{1.0ex plus .2ex}
\titlespacing*{\subsection}
	{0pt}{2.1ex plus .5ex minus .2ex}{.7ex plus .2ex}

\setlist{topsep=.45em,itemsep=.15em,parsep=0pt,partopsep=0pt}
\allowdisplaybreaks[2]
\numberwithin{equation}{section}

\newtheorem{theorem}{Theorem}[section]
\newtheorem{corollary}[theorem]{Corollary}
\newtheorem{lemma}[theorem]{Lemma}
\newtheorem{proposition}[theorem]{Proposition}

\newtheorem{definition}[theorem]{Definition}
\newtheorem{remark}[theorem]{Remark}
\newtheorem{example}[theorem]{Example}

\title{Codimension-Two Spacelike Submanifolds with Umbilical Lightlike Normal Sections and Their Relationship to Lightlike Hypersurfaces}

\author{
	Juan S. Gómez\thanks{\href{https://orcid.org/0000-0003-4038-2667}{ORCID: 0000-0003-4038-2667}} \\
	\href{mailto:jsalvador.gomez.edu@juntadeandalucia.es}{jsalvador.gomez.edu@juntadeandalucia.es} \\
	\small Departamento de Matemáticas \\
	\small Instituto E.S. Miguel Servet, 41020, Sevilla, Spain
}

\date{}
\begin{document}
	\maketitle
	\thispagestyle{empty}
	

	\begin{abstract}

We study codimension-two spacelike submanifolds in Lorentzian spacetimes that admit umbilical lightlike normal directions. We show that such submanifolds are subject to strong geometric and topological constraints, establishing explicit relationships between extrinsic geometry, mean curvature, and shear-isotropy. In the compact case, we obtain sharp restrictions on their topology. We precisely characterize when the umbilical lightlike normal vector field can be rescaled to be parallel, in terms of the curvature tensor of the ambient spacetime, and prove that this property is conformally invariant. Our main result is a factorization theorem: an embedded codimension-two spacelike submanifold carrying a lightlike normal vector field is contained, through the normal exponential construction, in a lightlike hypersurface endowed with a geodesic lightlike extension and an integrable screen distribution. The generated hypersurface is totally umbilical precisely when every leaf of this screen distribution is umbilical along the restriction of the lightlike generator. In particular, an umbilical lightlike normal section generates a totally umbilical lightlike hypersurface whenever the ambient spacetime is locally conformally flat. We also provide explicit conformal relations between the induced metrics on the family of spacelike leaves generated by the lightlike normal flow, with consequences for isometry, parallelism of lightlike normal directions, volume evolution, and variational properties. These results yield a detailed geometric framework relevant to the mathematical study of horizons and lightlike structures in general relativity.
	\end{abstract}
	
	\par
	
	\label{xyzt} 
	
\section{Introduction}

Lightlike hypersurfaces and their interplay with codimension-two spacelike submanifolds form a cornerstone of Lorentzian geometry, with far-reaching applications in general relativity, particularly in the study of black hole horizons~\cite{DuggalGimenez, Duggal_Sahin, kupeli}. In this paper, we examine how lightlike normal directions, especially umbilical ones, govern the extrinsic geometry of codimension-two spacelike submanifolds in a spacetime.

From a physical viewpoint, \emph{shear}, which measures the local instantaneous deformation, is intimately connected to horizons and their evolution. For example, isolated black hole horizons in equilibrium are Killing horizons~\cite{Wald}, characterized by vanishing shear. Mathematically, shear-free evolution corresponds precisely to the submanifold being umbilical along the evolution direction~\cite{cipriani}.

Motivated by these observations, the central concept in this paper is the \emph{umbilical lightlike normal section}. We prove that any codimension-two spacelike embedded submanifold carrying a lightlike normal vector field factors, via the normal exponential map, through a lightlike hypersurface endowed with a geodesic lightlike extension and an integrable screen distribution. The generated hypersurface is totally umbilical precisely when every leaf of this screen distribution is umbilical along the lightlike generator. In particular, an umbilical lightlike normal section generates a totally umbilical lightlike hypersurface whenever the ambient spacetime is locally conformally flat. This \emph{Factorization Theorem} (\Cref{factors}) provides an explicit link between the geometry of the initial spacelike leaf and that of the associated lightlike hypersurface, with applications to black hole horizons, cosmological models, and other lightlike structures in relativity~\cite{MinguzziSanchez}.

The paper is organized as follows. \Cref{sec:prelim} reviews background on Lorentzian geometry and submanifold theory. In \Cref{sec:shear_iso}, we study the relationship between umbilical lightlike normal sections and shear-isotropic spacelike submanifolds (see~\Cref{def:Shear_isotropy}), showing that, in the absence of umbilical points, shear-isotropy is equivalent to the existence of a globally defined umbilical lightlike normal section (\Cref{shear1theo} and~\Cref{shearCor}).

\Cref{sec:topo_ULS} explores extrinsic scalar curvature and mean curvature for codimension-two spacelike submanifolds with an umbilical lightlike normal section (\Cref{propo2}), and derives topological constraints on compact spacelike surfaces, which, under suitable hypotheses, must be tori or Klein bottles (Corollaries~\ref{GaussCurvature} and~\ref{TorusKlein}).

In \Cref{sec:para}, we give a curvature criterion for when an umbilical lightlike normal section can be locally rescaled to be parallel and show that, on a simply connected submanifold, the corresponding global curvature condition yields a global parallel rescaling and trivial normal holonomy~(\Cref{CriterionParaUmbi}~and~\Cref{PUS}). We also discuss examples where this condition holds and prove the conformal invariance of the rescalability property (\Cref{corollary:ParaConformal}).

Section~\ref{sec:hyper} develops the general theory of lightlike hypersurfaces needed for the Factorization Theorem: radical and screen distributions, the lightlike second fundamental form, and the totally umbilical case with its conformal invariance. These build up to the main result (\Cref{factors}):
\begin{quote}\itshape
	Let $\psi\colon S\to M$ be a codimension-two spacelike submanifold
embedded in a spacetime $M$, and let
$\xi\in\mathfrak X^\perp(S)$ be a lightlike normal vector field.
Then $\psi$ factors through a lightlike hypersurface $\Sigma$ of $M$
that admits a geodesic lightlike extension
$U\in\mathfrak X(\Sigma)$ of $\xi$ and an integrable screen
distribution $\mathcal S$ for which $S$ is one of the leaves.

If, moreover, $M$ is locally conformally flat and $\xi$ is an
umbilical section of $\psi$, then $\Sigma$ is totally umbilical.
\end{quote}

We also present a counterexample (\Cref{counterexample_GRW}) showing that even a stationary umbilical normal section may generate a totally umbilical lightlike hypersurface that is not totally geodesic, and conclude with remarks on the local and global applicability of the factorization framework, the submersion structure associated with the screen distribution, and the invariance properties of the resulting geometric structures.

Finally, \Cref{sec:conse_Umbil_Fac_Th} presents several remarkable consequences of the factorization under the explicit assumption that the generated lightlike hypersurface is totally umbilical. These include concrete conformal relations among the induced Riemannian metrics on leaves that project globally onto the initial submanifold~(\Cref{Cor:ConformalRelation}); invariance of the corresponding conformal factor under rescalings of the lightlike generator~(\Cref{corollary:invariance_omega_en}) and its precise transformation law under pointwise conformal changes of the ambient Lorentzian metric (see Remark~\ref{remark:Omega_conformal_invariance}); conditions for canonical isometry between the leaves~(Corollaries~\ref{Cor:isometry_condition} and~\ref{Cor:Omega_cte}); volume evolution formulas for compact leaves~(Theorems~\ref{Th:Vol1} and~\ref{propo:volume evolution formula}); variational and Jacobi field formulas~(\Cref{corSJacobi} and~\Cref{cor:umbilical-jacobi}); and conformally invariant curvature criteria for the local parallel rescalability of a lightlike vector field in the totally umbilical lightlike hypersurface~(\Cref{corollary:parallel_normal_connection}). When the ambient spacetime is locally conformally flat, the umbilicity of the initial lightlike normal section guarantees this standing hypothesis by~\Cref{factors}.

Together, these results provide a detailed geometric framework for umbilical lightlike normal sections and their associated totally umbilical lightlike hypersurfaces, clarifying the interplay between intrinsic and extrinsic properties of spacelike submanifolds in Lorentzian geometry, and offering tools of potential interest in rigorous studies of lightlike structures in general relativity.

	\section{Preliminaries}\label{sec:prelim}
	
	Let $(M,g)$ be a Lorentzian manifold with dimension ${\rm dim}\, M \geq 4$, that is,  $M$ is a connected smooth manifold and $g$ is a non-degenerate metric on $M$ with signature $(-,+,\ldots,+)$.  A tangent vector
	$v\in T_pM$, $p\in M$, is classified as spacelike (timelike or lightlike) if
		$g(v,v)>0$ or $v=0$ ($g(v,v)<0$, or $g(v,v)=0$ and $v\neq0$,
		respectively). A tangent vector is said to be causal if it is timelike or lightlike, that is, if $g(v,v)\leq 0$ and $v\neq 0.$

	A spacetime $(M,g)$ is defined as a time-orientable Lorentzian manifold endowed with a time orientation \cite{beem}. The manifold $M$ is not assumed to be  orientable in general. Note that orientability and time-orientability are logically independent.

		Throughout, all submanifolds are assumed to be connected unless otherwise stated.

		Let $\psi\colon S \to M$ be a codimension-two spacelike submanifold in an $(n+2)$-dimensional spacetime $(M,g).$ In other words, $S$ is an $n$-dimensional smooth manifold and $\psi$ is a smooth immersion such that the induced metric $\psi^*g$, also denoted by $g$, is Riemannian. For all local formulas and calculations, we may assume that $\psi$ is an embedding, and hence identify $p\in S$ with $\psi(p)\in M$, the tangent
	space $T_pS$ is identified with the subspace $\psi_*(T_pS)$ of $T_pM,$ and the normal space is denoted by $T_p^\perp S$. We will use letters $X,Y,Z$ (respectively $\xi,\eta,\zeta$) to represent vector fields tangent (respectively normal) to $S$.
	
	Throughout, we use the notation $(\cdot)^\top$ and $(\cdot)^\perp$ to denote the projections onto the tangent and normal bundles of $S$, respectively.
	
	Let $\widetilde{\nabla}$ and $\nabla$ be the Levi-Civita connections of $M$ and $S$, respectively. Then, the Gauss and Weingarten formulas are respectively expressed as follows:
	\begin{align} 
		\label{fgauss} 
		\widetilde{\nabla}_{X} Y &= \nabla_{X} Y + \mathrm{II}(X,Y),\\ 
		\widetilde{\nabla}_{X} \zeta &= - A_\zeta X + \nabla^\perp_X \zeta, \label{fWeingarten}
	\end{align} 
	where $\mathrm{II}$ represents the second fundamental form of $\psi,$  $A_\zeta$ the shape operator corresponding to $\zeta$, and $\nabla^\perp$ the normal connection. The shape operator and the second fundamental form are related by the following equation:
	\begin{equation} \label{relation}
		g( A_\zeta X, Y) = g( \mathrm{II}(X,Y), \zeta ).
	\end{equation}
	
	Denote the curvature tensors of $M$ and $S$ by $\widetilde{R}$ and $R$, respectively. The Gauss equation is given by
	\begin{equation}
		\bigl(\widetilde{R}(X,Y)Z\bigr)^{\top} = R(X,Y)Z + A_{{\rm II}(X,Z)}Y - A_{{\rm II}(Y,Z)}X, \label{GaussEc}
	\end{equation}
	for any $X, Y, Z \in {\mathfrak X}(S)$, where, according to our convention,
	\begin{equation*}
		R(X,Y)Z = \nabla_X \nabla_Y Z - \nabla_Y \nabla_X Z - \nabla_{[X,Y]}Z,
	\end{equation*}
	and similarly for  ${\widetilde R}$. Note that this is the opposite sign to the convention used by O'Neill~\cite[p.~74]{oneill}; we have chosen this sign for consistency with other references used in this work (e.g., \cite{beem, Dajczerbook, Duggal_Sahin}).
	
	Let $R^\perp$ denote the curvature tensor of the normal connection, that is,
	\[
	R^\perp(X, Y)\zeta = \nabla_X^\perp \nabla_Y^\perp \zeta - \nabla_Y^\perp \nabla_X^\perp \zeta - \nabla_{[X, Y]}^\perp \zeta
	\]
	for all $X, Y \in \mathfrak{X}(S)$ and $\zeta \in \mathfrak{X}^\perp(S)$. It follows from the Gauss and Weingarten formulas~\eqref{fgauss} and~\eqref{fWeingarten} that the normal component of $\widetilde{R}(X, Y)\zeta$ satisfies the Ricci equation:
	\begin{equation}\label{eq:RicciNormal}
		\bigl(\widetilde{R}(X, Y)\zeta\bigr)^\perp = R^\perp(X, Y)\zeta - \mathrm{II}(X, A_\zeta Y) + \mathrm{II}(Y, A_\zeta X).
	\end{equation}

	The \emph{mean curvature vector field of} $S$ in $M$ is defined by  $$\mathbf{H}=\frac{1}{n}\operatorname{trace}_g \,({\rm II})$$ 
	where $\operatorname{trace}_g({\rm II})$  is the $g$-trace of ${\rm II}$.

	Suppose that we can globally select two independent lightlike normal vector fields $\xi, \eta \in {\mathfrak X}^\perp(S)$ with $g(\xi, \eta) = -1$. Then the following (global) formulas hold:
	\begin{equation} \label{descomH}
		{\mathbf H} = -g(\eta, {\mathbf H})\, \xi - g(\xi, {\mathbf H})\, \eta,
	\end{equation}
	\begin{equation} \label{HExpansion}
		g({\mathbf H}, {\mathbf H}) = -2\,g(\xi, {\mathbf H})\, g(\eta, {\mathbf H}).
	\end{equation}
	Each component of ${\mathbf H}$ along the lightlike normal frame $\{\xi, \eta\}$ (up to a factor $n$) is referred to as the \emph{lightlike expansion scalar of $S$ along} $\xi$ and $\eta$, respectively:
	\begin{equation}\label{expansion-trace}
		\theta_\xi = -\operatorname{trace}(A_\xi) = -n\, g(\xi, {\mathbf H}), \qquad
		\theta_\eta = -\operatorname{trace}(A_\eta) = -n\, g(\eta, {\mathbf H}).
	\end{equation}
	
\begin{remark}\normalfont
	The term \emph{lightlike expansion scalar} originates from the physics literature~\cite{beem, HawkingEllis, kriele}, where it plays a central role in the analysis of lightlike congruences, particularly through the Raychaudhuri equation and its applications to the study of singularities and black hole horizons. The factor $1/n$ in the definition of the mean curvature vector is often omitted.
	
	The sign convention used in~\eqref{expansion-trace} for the mean curvature vector and the lightlike expansion scalar is determined by our choices in the Gauss and Weingarten formulas~\eqref{fgauss} and~\eqref{fWeingarten}~\cite{Dajczerbook}, and is consistent with the conventions adopted in~\cite{beem, DuggalGimenez, Duggal_Sahin, HawkingEllis}. This convention is used throughout the paper, including in the context of lightlike hypersurfaces discussed in later sections.
	
	We emphasize that this convention is standard in much of the mathematical literature, but it is the \emph{opposite sign} to that commonly used in the physics literature. All results and formulas in this paper are consistent with this choice.
	
	As sign conventions in the literature may vary depending on the adopted definitions, we hope that this explicit clarification will be helpful to readers from both the mathematical and physical communities.
\end{remark}

Whenever a global lightlike normal vector field $\xi$ exists, a lightlike normal frame $\{\xi, \eta\}$ with $g(\xi, \eta) = -1$ can always be constructed as follows. Take any globally defined timelike vector field $\mathbf{Z}$ on $M$ and, without loss of generality, assume that $\xi$ is future-directed, that is, $g(\xi, \mathbf{Z}) < 0$. Let $\mathbf{Z}^{\perp}$ denote the normal component of $\mathbf{Z}$ along $S$. Then define
\[
\eta = \frac{1}{g(\xi, \mathbf{Z})}\left( -\mathbf{Z}^{\perp} + \frac{g(\mathbf{Z}^{\perp}, \mathbf{Z}^{\perp})}{2g(\xi, \mathbf{Z})}\,\xi\right).
\]
This vector field satisfies $g(\eta, \eta) = 0$ and $g(\xi, \eta) = -1$, so $\{\xi, \eta\}$ forms a future-directed lightlike normal frame on $S$. Note that the resulting $\eta$ is unique and therefore independent of the particular choice of the global timelike vector field $\mathbf{Z}$.
	
	\begin{remark}\label{notaOrientable} \normalfont
		(a) The existence of a globally defined lightlike normal vector field on $S$ implies that its normal bundle is trivial. Moreover, if $M$ is orientable, then so is $S$~\cite[p.~214]{oneill}.
		
		(b) In general, even without a globally defined lightlike normal vector field, any codimension-two spacelike submanifold $S$ in a spacetime $M$ admits a globally defined, future-directed, timelike normal vector field: namely, the normal component $\mathbf{Z}^\perp$ of any global timelike vector field $\mathbf{Z}$ on $M$. In fact, $g(\mathbf{Z}^\perp, \mathbf{Z}^\perp) = g(\mathbf{Z}, \mathbf{Z}) - g(\mathbf{Z}^\top, \mathbf{Z}^\top) < 0$ along $S$, where $\mathbf{Z}^\top$ denotes the tangential component of $\mathbf{Z}$. However, the existence of a globally defined, non-vanishing spacelike normal vector field on $S$ is not always guaranteed; this is equivalent to the triviality of the normal bundle of $S$.
	\end{remark}

	Recall that a point $p \in S$ is said to be \emph{umbilical}~\cite{oneill} if there exists $\zeta \in T^{\perp}_p S$ such that
	\[
	{\rm II}(X, Y) = g(X, Y) \zeta,
	\]
	for any $X, Y \in T_p S$. The immersion $\psi$ is said to be \emph{totally umbilical} if every point of $S$ is umbilical. In this case, it is well known that
	\[
	{\rm II}(X, Y) = g(X, Y) \mathbf{H}
	\]
	for any $X, Y \in {\mathfrak X}(S)$.
	
	For a normal vector field $\zeta$ on $S$, if $A_\zeta = \rho I$ for some function $\rho$, then $\zeta$ is called an \emph{umbilical section} of $\psi$, or equivalently, $\psi$ is said to be umbilical with respect to $\zeta$. Hence, this function is given by $\rho = g(\zeta, \mathbf{H})$. Furthermore, from~\eqref{relation}, a submanifold is totally umbilical if and only if every normal section is umbilical. The spacelike submanifold $\psi$ is called \emph{pseudo-umbilical} if $\mathbf{H}$ is an umbilical section. Clearly, any totally umbilical spacelike submanifold is pseudo-umbilical.

	Observe that if $\psi$ is umbilical with respect to $\zeta$, then it is umbilical with respect to all vector fields proportional to $\zeta$. Therefore, we can say that $\psi$ is umbilical with respect to the normal direction spanned by $\zeta$, since this property concerns the umbilical direction $\mathrm{span}\{\zeta\}$, regardless of the length and whether $\zeta$ is future or past directed.
	
	\begin{remark} \label{CIS} \rm
		It is well known that the notion of umbilical section is invariant under conformal transformations of the ambient space~\cite{chen74}. In particular, it should be noted that the notion of an umbilical lightlike normal section is also conformally invariant.
	\end{remark}
	
	Consider $\zeta \in {\mathfrak X}^\perp(S)$ and the $n$-volume functional acting on compactly supported spacelike variations of $\psi$ along the normal direction defined by $\zeta$. The spacelike immersion $\psi$ is a critical point of this functional if and only if the component of the mean curvature vector field in the direction of $\zeta$ vanishes identically or, equivalently, if $g(\zeta, \mathbf{H}) = 0$~\cite[p.~299]{oneill}. Using equation~\eqref{expansion-trace}, this condition is equivalent to $\operatorname{trace}(A_\zeta) = 0$. In this case, the normal vector field $\zeta$ is referred to as \emph{stationary}.

	\section{Shear-isotropic Submanifolds}\label{sec:shear_iso}
	
	This section focuses on shear-isotropic submanifolds. Our analysis establishes a fundamental connection between shear-isotropy and the existence of lightlike umbilical sections, shedding light on the geometric features of these spacelike submanifolds.
	
	Let $\psi\colon S \to M$ be a codimension-two spacelike submanifold in a spacetime $M$. The \emph{total shear tensor} $\mathring{{\rm II}}$ is defined as the trace-free component of the second fundamental form~\cite{cipriani}:
	\begin{equation}\label{shear1}
		\mathring{{\rm II}}(X,Y) = {\rm II}(X,Y) - g(X,Y)\mathbf{H},
	\end{equation}
	for all $X, Y \in {\mathfrak X}(S)$. A point $p \in S$ is umbilical if and only if $\mathring{{\rm II}}_p = 0$.
	
	\begin{remark}\normalfont
		The total shear tensor is invariant under any conformal change of the ambient Lorentzian metric (see, e.g.,~\cite{fialkow}). Consequently, the set of umbilical points of $\psi$ is also a conformal invariant.
	\end{remark}
	
	The \emph{shear operator} associated with $\zeta \in {\mathfrak X}^\perp(S)$ is the trace-free part of the corresponding shape operator, that is,
	\begin{equation}\label{shear2}
		\mathring{A}_\zeta = A_\zeta - g(\zeta, \mathbf{H}) I.
	\end{equation}
	Using~\eqref{shear1} and~\eqref{shear2} in~\eqref{relation}, we see that the shear operators and the total shear tensor are related by the following equation:
	\begin{equation}\label{shear3}
		g(\mathring{A}_\zeta X, Y) = g(\mathring{\rm II}(X,Y), \zeta),
	\end{equation}
	for any $X, Y \in {\mathfrak X}(S)$.
	
	The \emph{shear scalar} $\sigma_\zeta$ associated with $\zeta \in {\mathfrak X}^\perp(S)$ is defined up to sign\footnote{The sign ambiguity is discussed in~\cite{cipriani, Senovilla}, but for our purposes only the value of $\sigma_\zeta^2$ is relevant.} as~\cite{Senovilla}
	\[
	\sigma_\zeta^2 = \operatorname{trace}(\mathring{A}_\zeta^2).
	\]
	Since $\mathring{A}_\zeta$ is self-adjoint, we have $\sigma_\zeta^2 \geq 0$, and $\sigma_\zeta^2 = 0$ if and only if $\mathring{A}_\zeta = 0$, or equivalently, $\zeta$ is an umbilical section of $\psi$.
	
	The \emph{shear space} of $\psi$ at $p \in S$ is defined as the subspace of $T_p^\perp S$ spanned by the values of the total shear tensor $\mathring{\rm II}$ at $p$, that is,
	\[
	\operatorname{Im}(\mathring{\rm II})_p = \operatorname{span}\left\{ \mathring{\rm II}(X,Y) : X, Y \in T_p S \right\}.
	\]
	A straightforward calculation shows that $\{\zeta \in T_p^\perp S : \mathring{A}_\zeta = 0\}$ is the orthogonal complement of $\mathrm{Im}(\mathring{\rm II})_p$ in $T_p^\perp S$.
	
	The concept of an isotropic submanifold, initially introduced by O'Neill in~\cite{oneill2} for Riemannian manifolds and subsequently generalized to pseudo-Riemannian manifolds in~\cite{Cab3}, is grounded in constraints on the second fundamental form. By replacing the latter with the total shear tensor~\eqref{shear1}, we can adapt the notion of isotropy and investigate the general geometric properties of these spacelike submanifolds.
	
	\begin{definition}\label{def:Shear_isotropy}\normalfont
		A codimension-two spacelike submanifold $\psi \colon S \to M$ in a spacetime $M$ is called \emph{shear-isotropic at a point} $p \in S$ if the real number
		\[
		g\bigl( \mathring{\rm II}(X,X), \mathring{\rm II}(X,X)\bigr) = \mathring{\lambda}(p)
		\]
		does not depend on the choice of unit tangent vector $X \in T_p S$. A submanifold is called \emph{shear-isotropic} if it is shear-isotropic at every point.
	\end{definition}
	
Clearly, this definition generalizes the notion of a totally umbilical spacelike submanifold.

Analogously to the case of isotropic spacelike surfaces in four-dimensional spacetime~\cite{Cab1}, we have the following result. 
	
	\begin{lemma} \label{CarShearIso}
		Let $\psi \colon S \to M$ be a codimension-two spacelike submanifold in a spacetime $M$. Then, the following conditions are equivalent:
		\begin{enumerate}
			\item $\psi$ is shear-isotropic.
			\item For any $X, Y, Z, W \in {\mathfrak X}(S)$, we have
			\[
			g\bigl(\mathring{{\rm II}}(X, Y), \mathring{{\rm II}}(Z, W)\bigr) = 0.
			\]
		\end{enumerate}
		
		Moreover, if $\psi$ is shear-isotropic at a non-umbilical point $p\in S$, then the shear space $\operatorname{Im}(\mathring{{\rm II}})_p$ is a  lightlike line of $T_p^\perp S$.
	\end{lemma}	
	\begin{proof}
		The proof is analogous to that of Theorem 4.2 in \cite{Cab1}. It is important to note that the value $\mathring{\lambda}$ in Definition~\ref{def:Shear_isotropy} is necessarily zero at each point of $S$, which is a crucial fact for the argument.
	\end{proof}
	
	It is clear from~\eqref{shear2} and~\eqref{shear3} that a codimension-two spacelike submanifold possessing an umbilical lightlike normal section is also shear-isotropic. Conversely, as a direct consequence of the following result, the shear-isotropy condition on a codimension-two spacelike submanifold guarantees the local existence of an umbilical lightlike section around any non-umbilical point. This intimate relationship between the two concepts provides valuable insight into the geometry of shear-isotropic submanifolds.
	
	\begin{theorem}\label{shear1theo}
		Let $\psi \colon  S \to M$ be a codimension-two spacelike submanifold in a spacetime $M$. Assume $\psi$ is free of umbilical points. If $\psi$ is shear-isotropic, then there exists a globally defined lightlike vector field $\xi \in {\mathfrak X}^\perp(S)$ and a smooth function $\rho \in C^\infty(S)$ such that $A_{\xi} = \rho I$. Furthermore, any other umbilical normal vector field is proportional to $\xi$ at each point of $S$. 
	\end{theorem} 
\begin{proof}
	Since $\psi$ is free of umbilical points, we can choose a locally finite open covering $\{\mathscr{U}_i\}$ of $S$ and, on each $\mathscr{U}_i$, a local $g$-orthonormal frame $E_{1}^i, \ldots, E_{n}^i \in \mathfrak{X}(\mathscr{U}_i)$ such that
	\[
	\xi_i = \mathring{\mathrm{II}}(E_{1}^i, E_{2}^i) \neq 0
	\]
	at every point of $\mathscr{U}_i$. By~\Cref{CarShearIso} and equation~\eqref{shear3}, we obtain
	\begin{align*}
		g(\xi_i, \xi_i) &= 0, \\
		g(\mathring{A}_{\xi_i}X, Y) &= 0,
	\end{align*}
	for all $X, Y \in \mathfrak{X}(S)$. Thus, $\xi_i$ is a lightlike normal section on $\mathscr{U}_i$ with $\mathring{A}_{\xi_i} = 0$.
	
	We may further arrange that $g(\mathbf{Z}, \xi_i) < 0$ on $\mathscr{U}_i$, where $\mathbf{Z}$ is a globally defined timelike vector field on $M$. Let $\{ f_i \}$ be a smooth partition of unity subordinate to the covering $\{\mathscr{U}_i\}$. Define $\xi = \sum_i f_i \xi_i$. The local finiteness of the sum ensures that $\xi$ is a well-defined smooth lightlike normal section on $S$. By~\Cref{CarShearIso}, at any $p \in \mathscr{U}_i \cap \mathscr{U}_j$ with $i \neq j$, we have $g(\xi_i, \xi_j) = 0$. Since both $\xi_i$ and $\xi_j$ are future-directed and lie in the same causal cone, which is convex, it follows that $\xi_i = a \xi_j$ for some constant $a > 0$. Consequently, at each $p \in S$, $\xi$ is non-zero because at least one function $f_i$ is strictly positive at $p$. Therefore, $\xi$ is a future-directed lightlike normal section on $S$ with $\mathring{A}_{\xi} = 0$.
	
	Finally, by~\Cref{CarShearIso} again, the shear space at $p \in S$ is the lightlike line
	\[
	\operatorname{Im}(\mathring{\mathrm{II}})_p = \left\{\zeta \in T_p^\perp S : \mathring{A}_\zeta = 0\right\}.
	\]
	Thus, if $\zeta$ is another normal vector at $p$ with $\mathring{A}_\zeta = 0$, then $\zeta \in \operatorname{Im}(\mathring{\mathrm{II}})_p = \operatorname{span}\{\xi\}$, so $\zeta$ and $\xi$ are collinear. This completes the proof.
\end{proof}

	The hypothesis that there are no umbilical points in~\Cref{shear1theo} is essential. Indeed, Example~5.5 of~\cite{Cab1} provides shear-isotropic spacelike surfaces in Lorentz-Minkowski spacetime $\mathbb{L}^4$ with umbilical points, where no global umbilical lightlike normal vector field exists.
	
	As a consequence of Theorem~\ref{shear1theo}, the following result holds:
	\begin{corollary}\label{shearCor}
		Let $S$ be a codimension-two spacelike submanifold in a spacetime $M$ without umbilical points. Then, $S$ is shear-isotropic and pseudo-umbilical if and only if there exists a stationary and umbilical lightlike normal section $\xi$ on $S$, that is, the associated shape operator satisfies $A_{\xi} = 0$.
	\end{corollary}
	\begin{proof}
		Suppose $S$ is shear-isotropic and pseudo-umbilical. By~\Cref{shear1theo}, there exists a lightlike normal vector field $\xi$ such that $A_{\xi} = \rho I$, where $\rho = g(\xi, \mathbf{H})$. Since $S$ is pseudo-umbilical, the mean curvature vector $\mathbf{H}$ is also umbilical. By the uniqueness of the umbilical direction (implied by~\Cref{shear1theo}), $\mathbf{H}$ must be collinear to $\xi$ at every point of $S$. Therefore, $g(\xi, \mathbf{H}) = 0$ and $A_{\xi} = 0$.
		
		Conversely, suppose that there exists an umbilical lightlike normal section $\xi$ on $S$ such that $A_\xi = 0$. Clearly, by~\eqref{shear3}, $S$ is shear-isotropic. Furthermore, since $\xi$ is stationary, $g(\xi, \mathbf{H}) = 0$. This implies that $\mathbf{H}$ is proportional to the umbilical section $\xi$. Hence, $S$ is pseudo-umbilical.
	\end{proof}

\begin{remark}\label{remark:stationary_umbilical}
	\normalfont
	(a) A codimension-two spacelike submanifold in a spacetime is isotropic if and only if it is both shear-isotropic and pseudo-umbilical, as can be verified using~\eqref{shear1} (see~\cite{ Cab3, Cab1, Cab2} for concrete examples and properties of isotropic submanifolds). Consequently, every isotropic codimension-two spacelike submanifold in a spacetime is inherently shear-isotropic.
	
	(b) If a codimension-two spacelike submanifold $S$ in a spacetime $M$ possesses a globally defined lightlike normal section $\xi$ on $S$ satisfying $A_{\xi} = 0$, equations~\eqref{HExpansion} and~\eqref{expansion-trace} imply that $g(\mathbf{H}, \mathbf{H}) = 0$. A particular case arises when $\mathbf{H}_p = 0$ at a point $p \in S$. In the specific scenario where $\mathbf{H}$ is lightlike everywhere, $S$ is called \emph{marginally trapped}~\cite{penrose}. Marginally trapped spacelike submanifolds are pivotal in general relativity, particularly in cosmology and black hole research (see, e.g., \cite{kupeli86}).
\end{remark}

	\section{Topology of Spacelike Surfaces with Umbilical Lightlike Sections}\label{sec:topo_ULS}
	
Next, we analyze the topological implications of umbilical lightlike sections on spacelike submanifolds in spacetimes. We show that the existence of such a section imposes stringent constraints on the topology of compact surfaces, restricting their possible shapes to tori or Klein bottles. These findings illustrate how local geometric conditions can influence the global topology of spacelike submanifolds.

Let \(\psi\colon S\to M\) be a codimension-two spacelike submanifold in a spacetime \(M\). For each point \(p\in S\) and each non-degenerate plane \(\Pi\subset T_pS\), we define the \emph{discriminant} at \(p\) along \(\Pi\) (cf.~\cite{oneill2}) by
\[
\mathscr{D}_p(\Pi)
= \mathscr{K}(\Pi) - \widetilde{\mathscr{K}}(\Pi),
\]
where \(\mathscr{K}(\Pi)\) and \(\widetilde{\mathscr{K}}(\Pi)\) are the sectional curvatures of \(S\) and \(M\), respectively.

The discriminant at $p \in S$ is said to be constant if $\mathscr{D}_p(\Pi)$ assumes the same value for all non-degenerate planes $\Pi \subset T_pS$. If this holds for every \(p\in S\), then $\mathscr{D}\colon S\to \mathbb{R}$ is a well-defined smooth function.

The \emph{extrinsic scalar curvature} of \(\psi\) is given by
\[
\tau_{\mathrm{ext}}
= \frac{2}{n(n-1)} \sum_{i<j} \mathscr{D}(\Pi_{ij}),
\]
where \(\Pi_{ij} = \mathrm{span}\{E_i, E_j\}\) for \(i \neq j\) and \(\{E_1, \dots, E_n\}\) is a local orthonormal frame on \(S\) (see~\cite{chen74} for a related definition in the Riemannian case). In the constant-discriminant case, \(\tau_{\mathrm{ext}} = \mathscr{D}\).

The next result expresses the extrinsic scalar curvature in terms of the mean curvature vector in the presence of an umbilical lightlike normal section.

\begin{proposition}\label{propo2}
	Let \(\psi\colon S\to M\) be a codimension-two spacelike submanifold in a spacetime \(M\) that admits an umbilical lightlike normal section. Then,
	\begin{equation}\label{ExtrinsicScalarH}
		\tau_{\mathrm{ext}} = g(\mathbf{H}, \mathbf{H}).
	\end{equation}
	
	Moreover, if \(M\) has constant sectional curvature \(c\), then the scalar curvature of $S$ satisfies
	\begin{equation}\label{Scal}
		\mathrm{Scal} = n(n-1)\bigl(c + g(\mathbf{H}, \mathbf{H})\bigr).
	\end{equation}
\end{proposition}
\begin{proof}
	From~\eqref{shear1}, we have
	\begin{align}
		g(\mathring{\mathrm{II}}, \mathring{\mathrm{II}})
		&= \sum_{i,j=1}^n g\bigl(\mathring{\mathrm{II}}(E_i, E_j),\, \mathring{\mathrm{II}}(E_i, E_j)\bigr) \nonumber \\
		&= \sum_{i,j=1}^n g\bigl(\mathrm{II}(E_i, E_j) - \delta_{ij} \mathbf{H},\, \mathrm{II}(E_i, E_j) - \delta_{ij} \mathbf{H}\bigr) \nonumber \\
		&= \sum_{i,j=1}^n g\bigl(\mathrm{II}(E_i, E_j),\, \mathrm{II}(E_i, E_j)\bigr) - 2\sum_{i=1}^n g\bigl(\mathrm{II}(E_i, E_i),\, \mathbf{H}\bigr) + n\, g(\mathbf{H}, \mathbf{H}) \nonumber \\
		&= g(\mathrm{II}, \mathrm{II}) - n\, g(\mathbf{H}, \mathbf{H}), \label{RelationSecond}
	\end{align}
	where $g(\mathrm{II}, \mathrm{II}) = \sum_{i,j=1}^n g\bigl(\mathrm{II}(E_i, E_j),\, \mathrm{II}(E_i, E_j)\bigr)$, and similarly for $g(\mathring{\mathrm{II}}, \mathring{\mathrm{II}})$. Both $g(\mathrm{II}, \mathrm{II})$ and $g(\mathring{\mathrm{II}}, \mathring{\mathrm{II}})$ are independent of the choice of local orthonormal frame $\{E_1, \ldots, E_n\}$ on $S$.
	
	Therefore, from the Gauss equation~\eqref{GaussEc} and equations~\eqref{relation} and~\eqref{RelationSecond}, we obtain
	\begin{align*}
		g(\mathring{\mathrm{II}}, \mathring{\mathrm{II}}) - n(n-1) g(\mathbf{H}, \mathbf{H}) 
		&= \sum_{i,j=1}^n g\bigl(\mathrm{II}(E_i, E_j), \mathrm{II}(E_i, E_j)\bigr) - \sum_{i,j=1}^n g\bigl(\mathrm{II}(E_i, E_i), \mathrm{II}(E_j, E_j)\bigr) \\
		&= \sum_{i \neq j} g\bigl(\mathrm{II}(E_i, E_j), \mathrm{II}(E_i, E_j)\bigr) - \sum_{i \neq j} g\bigl(\mathrm{II}(E_i, E_i), \mathrm{II}(E_j, E_j)\bigr) \\
		&= -n(n-1) \tau_{\mathrm{ext}}.
	\end{align*}
	This implies that
	\begin{equation} \label{RelationSecondBis}
		g(\mathring{\mathrm{II}}, \mathring{\mathrm{II}}) = n(n-1) \big( g(\mathbf{H}, \mathbf{H}) - \tau_{\mathrm{ext}} \big).
	\end{equation}
	
	Since $\psi$ admits an umbilical lightlike normal section $\xi$, we have $\mathring{A}_\xi = 0$. It then follows from~\eqref{shear3} that, at each point, the total shear tensor takes values in the lightlike direction spanned by $\xi$. Therefore, $g(\mathring{\mathrm{II}}, \mathring{\mathrm{II}}) = 0$, and from~\eqref{RelationSecondBis} we obtain~\eqref{ExtrinsicScalarH}.
	
	In particular, if $M$ is a spacetime of constant sectional curvature $c$, a straightforward computation yields~\eqref{Scal}.
\end{proof}
	
	We conclude this section by presenting a collection of results concerning spacelike surfaces in four-dimensional spacetimes.
	
	\begin{corollary}\label{GaussCurvature}
		Let $\psi \colon S \to M$ be a spacelike surface in a four-dimensional spacetime $M$ that admits an umbilical lightlike normal section. Then, the Gauss curvature of $S$ is given by
		\begin{equation} \label{curvature}
			\mathscr{K} = \widetilde{\mathscr{K}} + g(\mathbf{H}, \mathbf{H}),
		\end{equation}
		where $\widetilde{\mathscr{K}}$ denotes the sectional curvature of the ambient spacetime $M$ along $\psi$.
		
		Moreover, if $S$ is a compact spacelike surface, then 
		\begin{equation}\label{EulerChar}
			\chi(S) = \frac{1}{2\pi} \int_S \bigl( \widetilde{\mathscr{K}}+g(\mathbf{H}, \mathbf{H})  \bigr) \, {\rm d} \mu_g,
		\end{equation}
		where $\chi(S)$ is the Euler characteristic of $S$ and ${\rm d} \mu_g$ is the canonical measure associated with $(S,g)$.
	\end{corollary}
\begin{proof}
	By~\eqref{ExtrinsicScalarH}, it suffices to note that for a spacelike surface $S$ in a four-dimensional Lorentzian manifold $M$, the discriminant is constant, so $\tau_{\mathrm{ext}} = \mathscr{K} - \widetilde{\mathscr{K}}$ at each point of $S$. The formula~\eqref{EulerChar} then follows directly from~\eqref{curvature} by applying the Gauss--Bonnet Theorem.
\end{proof}
	
	\begin{remark}
		\normalfont
		Let $\psi \colon S \to M$ be a spacelike surface in a four-dimensional spacetime $M$ admitting an umbilical lightlike vector field $\xi \in {\mathfrak X}^\perp(S)$. Assume that $\psi$ is free of umbilical points. Then, there exists $\eta \in {\mathfrak X}^\perp(S)$ with $g(\eta,\eta)=0$ and $g(\xi, \eta)=-1$ such that the bilinear form $\mathring{\rm II}_\eta$, defined by
		\begin{equation*}
			\mathring{\rm II}_\eta(X,Y) = g(\mathring{A}_\eta X, Y)
		\end{equation*}
		for any $X,Y \in {\mathfrak X}(S)$, induces a Lorentzian metric on $S$. In fact, since ${\rm trace}(\mathring{A}_\eta) = 0$, the Cayley--Hamilton Theorem applied to the shear operator $\mathring{A}_\eta$,
		\[
		\mathring{A}_\eta^2 - \operatorname{trace}(\mathring{A}_\eta) \mathring{A}_\eta + \operatorname{det}(\mathring{A}_\eta)I = 0,
		\]
		yields
		\[
		2 \operatorname{det}(\mathring{A}_\eta) = - \operatorname{trace}(\mathring{A}_\eta^2) \leq 0,
		\]
		with equality if and only if $\eta$ is an umbilical lightlike section of $\psi$.
	\end{remark}

	Recall that a connected manifold admits a Lorentzian metric if and only if it admits a non-vanishing tangent vector field~\cite[Prop.~5.37]{oneill}. By elementary algebraic topology, this property is equivalent to $S$ being either non-compact, or compact with the Euler characteristic $\chi(S) = 0$. Therefore, the only compact surfaces that can admit Lorentzian metrics are the torus and the Klein bottle.
	
	As a direct consequence of Corollary~\ref{GaussCurvature} and the classification of compact surfaces admitting Lorentzian metrics, we obtain the following:
	
	\begin{corollary}\label{TorusKlein}
		Let $\psi \colon S \to M$ be a compact spacelike surface in a four-dimensional spacetime $M$ that admits an umbilical lightlike normal section and is free of umbilical points. Then
		\begin{equation*}
			\int_S \bigl(\widetilde{\mathscr{K}} + g(\mathbf{H}, \mathbf{H}) \bigr)\,{\rm d} \mu_g = 0,
		\end{equation*}
		and $S$ is diffeomorphic to either the torus or the Klein bottle.
	\end{corollary}

	This result reveals a profound connection between the intrinsic and extrinsic geometry of compact spacelike surfaces in four-dimensional spacetimes. In particular, it imposes strong topological constraints, restricting the possible diffeomorphism types of such surfaces to tori or Klein bottles. The implications of this finding extend to various areas of theoretical physics, including cosmology and string theory (see, e.g.,~\cite{GonzalezAlonso}).
	
		On the other hand, we have the following result for the $2$-sphere:

\begin{corollary}\label{UmbilicalPoint}
	Every spacelike surface in a four-dimensional spacetime, which is topologically a $2$-sphere and admits an umbilical lightlike normal section, possesses at least one umbilical point.
\end{corollary}
	
	Every compact spacelike surface $S$ in Lorentz-Minkowski spacetime $\mathbb{L}^4$ that factors through a light cone possesses a lightlike umbilical normal section and is diffeomorphic to a $2$-sphere~\cite{PalomoRomero13}. Consequently,~\Cref{UmbilicalPoint} guarantees the existence of at least one umbilical point $p \in S$ (cf.~Remark~2.3 in~\cite{PalomoRomero14}).

	\section{Parallelism of Umbilical Lightlike Normal Sections}\label{sec:para}
	The aim of this section is to investigate the conditions under which an umbilical lightlike normal section is parallel. Such an analysis is fundamental for understanding the relationship between these spacelike submanifolds and lightlike hypersurfaces in a spacetime.
	
	Let $\psi\colon S \to M$ be a codimension-two spacelike submanifold in a spacetime $M$. Suppose that $\{\xi, \eta\}$ is a globally defined normal lightlike frame on $S$ such that $g(\xi, \eta) = -1$. Since $g(\xi, \xi) = 0$, it follows that $g(\nabla_X^\perp \xi, \xi) = 0$ for every $X \in \mathfrak{X}(S)$. Consequently, there exists a one-form $\tau$ on $S$ such that $\nabla_X^\perp \xi = \tau(X)\xi,$ given by  
	\begin{equation}\label{eq:tau_definition}
		\tau(X) = -g(\widetilde{\nabla}_X \xi, \eta).
	\end{equation}
	
	Using the well-known Maurer-Cartan formula \[ d\tau(X, Y) = X(\tau(Y)) - Y(\tau(X)) - \tau([X, Y]), \]
	it is straightforward to verify that
	\[
	R^{\perp}(X, Y) \xi = d\tau(X, Y)\, \xi, \quad 
	R^{\perp}(X, Y) \eta = -d\tau(X, Y)\, \eta,
	\]
	for any $X, Y \in \mathfrak{X}(S)$. Therefore, for every normal vector field $\zeta \in \mathfrak{X}^\perp(S)$, we obtain 
	\begin{equation}\label{eq:RNorTau}
		R^\perp(X,Y)\zeta = d\tau(X,Y)\overline{\zeta},
	\end{equation}
	where $\overline{\zeta} = -g(\zeta, \eta)\xi + g(\zeta,\xi)\eta$.

	\begin{remark}
		\normalfont
		If a codimension-two spacelike submanifold $S$ in a spacetime $M$ admits a non-vanishing umbilical normal vector field, then the shape operators associated with any two normal vector fields commute at every point of $S$. This is a direct consequence of the basic properties of the shape operators. In this case, taking into account~\eqref{relation}, the Ricci equation~\eqref{eq:RicciNormal} simplifies to
		\begin{equation}\label{eq:RSpacetimePerp}
			\bigl(\widetilde{R}(X,Y)\zeta\bigr)^\perp = R^\perp(X,Y)\zeta,
		\end{equation}
		for any $X, Y \in \mathfrak{X}(S)$ and $\zeta \in \mathfrak{X}^\perp(S)$.
	\end{remark}

	\begin{proposition}\label{CriterionParaUmbi}
		Let $\psi\colon S \to M$ be a codimension-two spacelike submanifold in a spacetime $M$. The following conditions are equivalent:
		\begin{enumerate}
			\item There exists a globally defined umbilical lightlike normal vector field $\xi$ on $S$, that can be locally rescaled to be parallel with respect to the normal connection.
			\item There exists a globally defined umbilical lightlike normal vector field $\xi$ on $S$ such that
			\begin{equation}\label{CurvatureRestric}
			\bigl(\widetilde{R}(X,Y)\xi\bigr)^\perp = 0,
		\end{equation}
			for all $X, Y \in \mathfrak{X}(S)$, where $\widetilde{R}$ is the curvature tensor of $M$.
		\end{enumerate}
		In this case, the submanifold $S$ has a flat normal connection.
	\end{proposition}
	
	\begin{proof}
		Assume first that $\xi$ is an umbilical lightlike normal vector field satisfying~\eqref{CurvatureRestric}. Then, by~\eqref{eq:RNorTau} and~\eqref{eq:RSpacetimePerp}, the one-form $\tau$ is closed ($d\tau = 0$), and thus, by the Poincaré Lemma, for each point $p \in S$ there exists an open neighborhood $\mathscr{U}$ of $p$ and a smooth function $f$ on $\mathscr{U}$ such that $\tau = df$. If $S$ is simply connected, then $\tau = df$ globally on $S$. In this case, by defining $\overline{\xi} = e^{-f} \xi$, a globally defined parallel umbilical lightlike normal section $\overline{\xi}$ can be constructed, as we show for the local case below.
		
		Fix a point $p \in S$. Let $\mathscr{U}$ be an open neighborhood of $p$ where there exists a smooth function $f$ such that $\tau = df$. Define the local lightlike normal vector fields $\overline{\xi} = e^{-f} \xi$ and $\overline{\eta} = e^{f} \eta$ on $\mathscr{U}$. Note that $g(\overline{\xi}, \overline{\eta}) = -1$, and for any $X \in \mathfrak{X}(\mathscr{U})$,
		\begin{align*}
			\overline{\tau}(X) &= -g(\widetilde{\nabla}_X \overline{\xi}, \overline{\eta}) \\
			&= -g\bigl(\widetilde{\nabla}_X (e^{-f} \xi), e^{f} \eta\bigr) \\
			&= X(f)g(\xi,\eta) - g(\widetilde{\nabla}_X \xi, \eta) \\
			&= -X(f) + \tau(X) \\
			&= 0,
		\end{align*}
		since $\tau(X) = df(X) = X(f)$. Therefore, $\overline{\xi}$ is a parallel umbilical lightlike normal vector field on $\mathscr{U}$.
		
		Conversely, suppose that there exists a globally defined umbilical lightlike normal vector field $\xi$ on $S$ that can be locally rescaled to be parallel with respect to the normal connection. That is, for each point $p \in S$, there exists an open neighborhood $\mathscr{U}$ of $p$ and a smooth function $f$ on $\mathscr{U}$ such that the rescaled field $\overline{\xi} = e^{-f}\xi$ is parallel on $\mathscr{U}$. Since the associated one-form $\overline{\tau}$ vanishes on $\mathscr{U}$, and given that $\overline{\tau} = \tau - df$ as shown above, it follows that $\tau = df$ on each $\mathscr{U}$. Therefore, $d\tau = 0$ on $S$. Substituting $d\tau = 0$ into~\eqref{eq:RNorTau} and~\eqref{eq:RSpacetimePerp} yields~\eqref{CurvatureRestric}.
		
		Finally, the fact that $d\tau = 0$ shows, by~\eqref{eq:RNorTau}, that the normal connection of $S$ is flat.
	\end{proof}

	Proposition~\ref{CriterionParaUmbi} gives the local curvature criterion for parallel rescalability. On a simply connected submanifold, the corresponding rescaling globalizes.

\begin{theorem}\label{PUS}
	Let $S$ be a codimension-two spacelike submanifold in a spacetime $M$ that admits an umbilical lightlike normal section $\xi$. Suppose that $S$ is simply connected. If $\bigl(\widetilde R(X,Y)\xi\bigr)^\perp=0$ for all $X,Y\in\mathfrak X(S)$, then $\xi$ can be rescaled to be parallel with respect to the normal connection on $S$. Moreover, the submanifold $S$ has a flat normal connection and the normal holonomy group of $S$ is trivial.
\end{theorem}

\begin{proof}
	By~\Cref{CriterionParaUmbi}, the curvature assumption implies $d\tau=0$, where $\nabla_X^\perp\xi=\tau(X)\xi$, and the normal connection of $S$ is flat. Since $S$ is simply connected, there exists $f\in C^\infty(S)$ such that $\tau=df$. Therefore, $\widehat\xi=e^{-f}\xi$ is a globally defined parallel lightlike normal vector field. Finally, a flat connection over a simply connected manifold has trivial holonomy.
\end{proof}
\begin{example}
	\normalfont
	In a locally conformally flat spacetime $M$, any codimension-two spacelike submanifold $S$ admitting an umbilical lightlike normal section possesses, \emph{locally}, a parallel umbilical lightlike normal field after a suitable rescaling. Indeed, the conformally flat curvature decomposition and the Ricci equation~\eqref{eq:RSpacetimePerp} show that the normal connection is flat in the presence of the umbilical normal direction; the conclusion then follows from~\Cref{CriterionParaUmbi}. If $S$ is simply connected, the rescaling can be performed globally by~\Cref{PUS}. Examples of locally conformally flat spacetimes include Lorentz--Minkowski, de Sitter, anti-de Sitter, and Robertson--Walker spacetimes~\cite{HawkingEllis,Wald}.
\end{example}

\begin{remark}\normalfont
	There exist spacetimes that are not locally conformally flat but nevertheless admit codimension-two spacelike submanifolds endowed with a parallel umbilical lightlike normal section. A notable example is provided by generalized Schwarzschild spacetimes~\cite{moron-palomo24}. These spacetimes are of particular physical and geometric interest, as they generalize black hole models to settings with more intricate spatial geometry and causal structure.
\end{remark}

Consider a metric \( g^* = e^{2u} g \) conformally equivalent to \( g \), where \( u \) is a smooth function on \( M \). Analogously to the Riemannian context~\cite{chen74}, we note that while the parallelism of a normal vector field is not, in general, preserved under conformal transformations, the normal curvature tensor \( R^\perp \) does possess conformal invariance. This is consistent with the conformal nature of the notion of an umbilical lightlike normal section (see~\Cref{CIS}). Therefore, by~\Cref{CriterionParaUmbi}, together with~\eqref{eq:RNorTau} and~\eqref{eq:RSpacetimePerp}, we obtain the following significant consequence:

\begin{corollary}\label{corollary:ParaConformal}
	The property that an umbilical lightlike normal vector field on a codimension-two spacelike submanifold in a spacetime can be rescaled to be parallel is invariant under conformal transformations of the ambient spacetime.
\end{corollary}



	\section{Umbilical Lightlike Normal Sections and Totally Umbilical Lightlike Hypersurfaces}\label{sec:hyper}

	We now explore the relationship between codimension-two spacelike submanifolds and lightlike hypersurfaces. We show that any codimension-two spacelike submanifold embedded in a spacetime that possesses a lightlike normal vector field can be embedded in a lightlike hypersurface. Although this result is known (see \cite[p.~31]{kupeli}), we provide an explicit construction, offering crucial insights into the geometry of such submanifolds and their connection to lightlike hypersurfaces. In what follows, we adopt the terminology of~\cite{Bejancu1996}.
	
	\subsection{Geometry of Lightlike Hypersurfaces}
	
	A \emph{lightlike hypersurface} in a spacetime $M$ is a codimension-one submanifold $\Sigma$ such that the induced metric on $\Sigma$ is everywhere degenerate. The \emph{radical} of $\Sigma$ at a point $p \in \Sigma$ is defined as $\mathrm{Rad}_p(\Sigma) = T_p\Sigma \cap T_p^\perp \Sigma,$ where $T_p\Sigma$ is the tangent space to $\Sigma$ at $p$ and $T_p^\perp \Sigma$ is its normal space \cite[p.~53]{oneill}. Since $T_p\Sigma$ is a degenerate hyperplane and the ambient metric is Lorentzian, we have $T_p^\perp \Sigma \subset T_p\Sigma$ and $\dim T_p^\perp \Sigma = 1$. Consequently, there exists a unique lightlike direction in $\Sigma$ that is orthogonal to any direction in $\Sigma$. In particular, $\Sigma$ does not contain any timelike directions and is foliated by a family of lightlike curves.

	\begin{remark} \normalfont \label{Foliation}
		The foliation of a hypersurface by lightlike curves is a necessary, but not sufficient, condition for it to be lightlike. For example, a timelike hyperplane in Lorentz-Minkowski spacetime is foliated by lightlike curves, yet it is not a lightlike hypersurface. However, if we impose an additional causality condition, namely local achronality, then the foliation by lightlike curves \emph{does} guarantee that the hypersurface is lightlike. Specifically, a locally achronal hypersurface in a spacetime that is foliated by lightlike curves is necessarily a lightlike hypersurface (see~\cite{kupeli}).
	\end{remark}
	
Every lightlike hypersurface $\Sigma$ in a spacetime $M$ gives rise to a future-directed lightlike vector field $U \in \mathfrak{X}(\Sigma)$ that generates the radical distribution $\mathrm{Rad}(\Sigma)$ and satisfies $\widetilde{\nabla}_{U} U = fU$ for some smooth function $f\colon \Sigma \to \mathbb{R}$~\cite[p.~116]{DuggalGimenez}; see also~\cite[Sec.~3]{kupeli87}. Furthermore, the lightlike vector field $U$ is unique up to a positive pointwise scale factor.
	
	An inextendible integral curve $\gamma\colon I \to \Sigma$ of $U \in \mathfrak{X}(\Sigma)$ (up to parametrization) is called a \emph{lightlike generator} of $\Sigma$ \cite{kupeli87}. 
	Note that every lightlike generator of $\Sigma$ is a lightlike pregeodesic in $M$, that is, it can be reparameterized to be a lightlike geodesic in $M$ (see, e.g.,~\cite{ MinguzziSanchez}). 
	
	The lightlike vector field $U$ is said to be \emph{geodesic} if it satisfies $\widetilde{\nabla}_{U} U = 0$. In this case, the lightlike generators are geodesics of ambient spacetime $M$. It can always be done locally, but, in general, it may not be possible to rescale the vector field $U$ to be a geodesic vector field.

	\begin{remark}\normalfont \label{remark:causally_separated}
		A classical question in the geometry of lightlike hypersurfaces concerns the global rescaling of a lightlike vector field to a geodesic one, a problem addressed in \cite[Sec.~4]{kupeli87}. A key enabling condition for this is the existence of a spacelike hypersurface $S \subset \Sigma$ that intersects each lightlike generator of $\Sigma$ at exactly one parameter value, ensuring the rescalability of any lightlike vector field $U \in \mathfrak{X}(\Sigma)$ to be geodesic. In this case, $\Sigma$ is diffeomorphic to the product manifold $S \times \mathbb{R}$.
	\end{remark}

	The \emph{lightlike second
		fundamental form} of $\Sigma$ associated with $U$ is the tensor field defined by \cite{Bejancu1996} 
	\begin{equation}\label{DefLightlikeSecond}
		B_U(X,Y)=-g(\widetilde{\nabla}_XU,Y)
	\end{equation}
	for any $X,Y \in {\mathfrak X}(\Sigma)$. Using the fact that $[X,Y] \in {\mathfrak X}(\Sigma)$ for any $X,Y \in {\mathfrak X}(\Sigma)$, it is straightforward to verify that $B_U$ is symmetric and $B_U(X,U) = 0$ for any $X \in {\mathfrak X}(\Sigma)$. The minus sign in \eqref{DefLightlikeSecond} is chosen for consistency with our convention for spacelike submanifolds (see~\eqref{fgauss}--\eqref{relation}).
	
	To  handle the degenerate metric of a lightlike hypersurface, it is useful to introduce a screen distribution: a non-degenerate spacelike subbundle of the tangent bundle that, together with the radical, spans the entire tangent space. This decomposition allows us to separate the problematic lightlike direction from the well-behaved spacelike directions, facilitating the analysis of induced geometric objects.
	
	Formally, a \emph{screen distribution} $\mathcal{S}$ is a complementary distribution to $\mathrm {Rad}(\Sigma)$ in $T\Sigma$. Then, $\mathcal{S}$ is necessarily spacelike distribution and we have orthogonal direct sum 
	\begin{equation}\label{DesRad}
		T\Sigma= \mathcal{S}\oplus {\rm Rad}(\Sigma).
	\end{equation}
	As $\Sigma$ is assumed to be paracompact, there is always ${\mathcal S}$ \cite{Duggal_Sahin}. 
	
	The \emph{transverse distribution} is the unique one-dimensional lightlike distribution orthogonal to $\mathcal{S}$ and not contained in $T\Sigma$ \cite[p.~44]{Duggal_Sahin}. Since $M$ is time-oriented, there exists a lightlike vector field $V$ over $\Sigma$ that generates the transverse distribution and can be normalized so that $g(U, V)=-1.$ This vector field is called the \emph{lightlike transverse vector field} associated with $U$ (and corresponding to $\mathcal{S}$).

	\begin{remark} \normalfont \label{remark:canonical_screen_Kupeli}
		Note that, in general, the screen distribution $\mathcal{S}$ is not canonical (thus not unique), and the lightlike geometry depends on its choice, an issue that has been extensively studied in the literature (see, e.g., \cite{DuggalGimenez} and \cite{Duggal_Sahin} and references therein). Although Kupeli~\cite{kupeli87} provided a canonical perspective and showed that $\mathcal{S}$ is isometric to the quotient vector bundle $T\Sigma/\mathrm{Rad}(\Sigma)$ for intrinsic geometric studies, our investigation adopts an extrinsic viewpoint as in Bejancu~\cite{Bejancu1996}, which is in line with the classical theory of non-degenerate submanifolds \cite{Dajczerbook}. Consequently, for the lightlike hypersurface constructed in relation to the codimension-two spacelike submanifold, an associated screen distribution with notable geometric properties inherently exists, a crucial element in exploring their relationship, as we shall see.
	\end{remark}
	
	We denote by $\Gamma({\mathcal S})$ the set of smooth sections of a screen distribution $\mathcal{S}$. Given $X\in {\mathfrak X}(\Sigma)$, the vector field $\widetilde{\nabla}_X U$ belongs to $\mathfrak{ X}(\Sigma)$ (see, e.g., \cite{Bejancu1996}), so it can be decomposed from \eqref{DesRad} as
	\begin{equation}\label{LightlikeDes1}
		\widetilde{\nabla}_X U= -A^*_UX+\tau(X)U,
	\end{equation}
	where $A^*_UX \in \Gamma({\cal S})$ and $\tau$ is the \emph{rotation one-form} given by
	\begin{equation}\label{rotation}
		\tau(X)=-g(\widetilde{\nabla}_X U,V).
	\end{equation}
	The endomorphism $A^*_U$ is called the \emph{screen shape operator} of the screen distribution $\mathcal S$ and, from \eqref{DefLightlikeSecond} and \eqref{LightlikeDes1}, satisfies 
	\begin{equation}\label{RelaLightlik}
		B_U(X,Y)=g(A^*_UX,Y)
	\end{equation}
	for any $X,Y \in {\mathfrak X}(\Sigma)$. 
	
\begin{remark}\label{remark:dtau_ditinghished}\normalfont 
	The rotation one-form $\tau$, which is determined by the chosen screen distribution, plays an important role in the geometry of lightlike hypersurfaces \cite{DuggalSharma}, particularly when it is closed, that is, $d\tau=0$. Observe that while $\tau$ depends on the choice of $U$, its exterior derivative $d\tau$ remains invariant under such choices~\cite[Prop.~2.3.1]{Duggal_Sahin}.
\end{remark}
	
	 A screen distribution ${\mathcal S}$ is said to be \emph{$U$-distinguished} if $\widetilde{\nabla}_X U \in \Gamma({\mathcal S})$ for any section $X$ of ${\mathcal S}$ \cite{DuggalGimenez}, a condition equivalent to $\tau(X)=0$ for all $X\in\Gamma({\mathcal S})$, as follows from~\eqref{LightlikeDes1}. Furthermore, if $U$ is also geodesic, then the associated one-form $\tau$ vanishes identically \cite[Prop.~4]{DuggalGimenez}.

Consider the restriction of the screen shape operator
$
A^{*}_{U}\colon \Gamma(\mathcal{S}) \to \Gamma(\mathcal{S}).
$
This operator is self-adjoint and hence diagonalizable, and its eigenvalues are independent of the chosen screen distribution (see,~e.g., \cite[Lem.~1]{DuggalGimenez}). It follows that the \emph{lightlike expansion scalar} of \(\Sigma\) with respect to \(U\), defined by
$
\theta_{U}= -\operatorname{trace}(A^{*}_{U}),
$
is likewise independent of that choice.
	
	\begin{remark} \normalfont \label{remark:sign_theta}
		Note that under a rescaling $\overline{U} = fU$ by a non-vanishing smooth function $f$ on $\Sigma$, the corresponding lightlike expansion scalar transforms as $\theta_{\overline{U}}= f\theta_U$. Consequently, its sign remains invariant under positive rescalings of the future-directed lightlike vector field $U$. Furthermore, the property of having a vanishing lightlike expansion scalar is independent of the chosen lightlike section $U$ on the hypersurface. Therefore, it is meaningful to refer to \emph{expansion-free} lightlike hypersurfaces as those possessing a vanishing lightlike expansion scalar.
	\end{remark}
	
	\begin{example} \normalfont In Lorentz-Minkowski spacetime endowed with the usual time-orientation \cite{oneill}, with the sign convention adopted herein for the lightlike expansion scalar \cite{DuggalGimenez}, a future light cone exhibits positive lightlike expansion scalar (relative to its future-directed lightlike generators), while a past light cone exhibits negative lightlike expansion scalar (see \cite{Bejancu1996} for details).
	\end{example}
	
The \emph{shear scalar} $\sigma_U$ associated with $U$ is defined (up to sign) as $\sigma_U^2 = \operatorname{trace}(\mathring{A}^*_U \circ \mathring{A}^*_U)$, where $\mathring{A}^*_U\colon \Gamma(\mathcal{S}) \to \Gamma(\mathcal{S})$ denotes the trace-free part of the screen shape operator, that is,
\begin{equation}\label{LightlikeShear1}
	\mathring{A}^*_U = A^*_U + \frac{\theta_U}{n} I,
\end{equation}
where $\theta_U$ is the lightlike expansion scalar defined above. Given that the eigenvalues of the screen shape operator $A^*_U$ are independent of the choice of screen distribution, any symmetric function of these eigenvalues is also independent of this choice. Expanding $A^*_U \circ A^*_U = \bigl(\mathring{A}^*_U - \frac{\theta_U}{n} I\bigr)^2$ and taking the trace, and noting that $\operatorname{trace}(\mathring{A}^*_U) = 0$, we obtain
\[
\operatorname{trace}(A^*_U \circ A^*_U) = \operatorname{trace}(\mathring{A}^*_U \circ \mathring{A}^*_U) + \frac{\theta_U^2}{n} = \sigma_U^2 + \frac{\theta_U^2}{n}.
\]
Therefore, the shear scalar $\sigma_U$ is also independent of the choice of screen distribution.  Furthermore, the property of having vanishing shear scalar is independent of the choice of lightlike section $U$, since under a rescaling $\overline{U} = fU$ by a non-vanishing smooth function $f$, the shear scalar transforms as $\sigma^2_{\overline{U}} = f^2 \sigma^2_U$. Thus, it makes sense to refer to \emph{shear-free} lightlike hypersurfaces.
	
	\subsection{Totally Umbilical Lightlike Hypersurfaces}
	\begin{definition}\normalfont \cite{Bejancu1996}
		A lightlike hypersurface $\Sigma$ in $M$ is called \emph{totally umbilical} if there exists a smooth function $\mu$ on $\Sigma$ such that 
		\begin{equation} \label{TUL}
			B_U(X,Y)=\mu g(X,Y),
		\end{equation}
		for any $X, Y\in {\mathfrak X}(\Sigma)$. Hence, this function $\mu$ is given by $\mu=-\theta_U/n$. In particular, $\Sigma$ is called \emph{totally geodesic} if $B_U$ vanishes.
	\end{definition}
	
	Note that the above definition does not depend on the particular choice of $U$. It is important to mention that the second fundamental form $B_U$ on $\Sigma$ is independent of the
	choice of screen distribution \cite[Prop.~2.1]{Bejancu1996}. Therefore, $\Sigma$ is totally umbilical if and only if \eqref{TUL} holds for all sections of a given screen distribution.
	
	\begin{remark}\normalfont
		Since the operator defined by equation \eqref{LightlikeShear1} is self-adjoint and diagonalizable, it follows from \eqref{RelaLightlik} that a lightlike hypersurface is totally umbilical if and only if it is shear-free.  Likewise, it is totally geodesic if and only if it is expansion- and shear-free.
	\end{remark}
	
\begin{example}
	\normalfont
	Lightlike hyperplanes in Lorentz-Minkowski spacetime are totally geodesic, whereas light cones are totally umbilical lightlike hypersurfaces that are not totally geodesic~\cite{Bejancu1996}. Event horizons of stationary black holes, including those of Kruskal, Reissner-Nordström, and Kerr, as well as compact Cauchy horizons in spacetimes satisfying the weak energy condition (such as the Taub-NUT spacetime), are examples of totally geodesic lightlike hypersurfaces~\cite[p.~323]{HawkingEllis}. Totally umbilical lightlike hypersurfaces in generalized Robertson--Walker spacetimes were studied and classified in~\cite{gutierrez15}.
\end{example}

\begin{remark}\label{remark:TU_Hyper} \normalfont
		Totally geodesic lightlike hypersurfaces have a geometric interpretation similar to the Riemannian case. Specifically, if $\gamma$ is a geodesic of $M$ starting tangent to $\Sigma$, then $\gamma$ remains in $\Sigma$ for sufficiently small parameter values. This follows from the fact that, when $\Sigma$ is totally geodesic, the restriction to $\Sigma$ of the Levi-Civita connection of $M$ defines an affine connection on $\Sigma$ (see~\cite{kupeli87} for details).
	
	In contrast, the totally umbilical hypothesis is different because the lightlike second fundamental form has a distinguished direction that belongs to its kernel.
	
In general, for a totally umbilical lightlike submanifold $\Sigma$ of a Lorentzian manifold $M$, it holds that, for some parameter interval around the starting point, a lightlike geodesic of $M$ that starts tangent to $\Sigma$ remains in $\Sigma$ (see, e.g.,~\cite[Prop.~3.4.3]{Duggal_Sahin}). This condition does not imply that $\Sigma$ is totally umbilical \cite{Perlick}. For lightlike hypersurfaces in a spacetime, however, this property is automatic for any lightlike geodesic starting tangent to the hypersurface, regardless of whether it is totally umbilical or not, because the lightlike generators are always pregeodesics of the ambient spacetime.
\end{remark}

	The following well-known result is straightforward to prove. For completeness, we include a proof. 
	
	\begin{lemma}\label{lemmaULH} 
		Let $\Sigma$ be a lightlike hypersurface in a spacetime $M$ and let $U \in {\mathfrak X}(\Sigma)$ be a lightlike vector field on $\Sigma$. The following conditions are equivalent:
		\begin{enumerate}
			\item $\Sigma$ is totally umbilical (respectively, totally geodesic).
			\item $U$ is a conformal (respectively, Killing) vector field of the degenerate induced metric on $\Sigma$.
		\end{enumerate}
	\end{lemma}
\begin{proof}
	Using the expression for the Lie derivative and the compatibility of the metric, we have
	\begin{align*}
		\mathscr{L}_U g(X,Y) &= U\bigl(g(X,Y)\bigr) - g(\mathscr{L}_U X, Y) - g(X, \mathscr{L}_U Y) \\
		&= g(\widetilde{\nabla}_X U, Y) + g(X, \widetilde{\nabla}_Y U)
	\end{align*}
	for any $X, Y \in {\mathfrak X}(\Sigma)$. The symmetry of the lightlike second fundamental form $B_U$ of $\Sigma$ associated with $U$ then gives $
	\mathscr{L}_U g(X,Y) = -2 B_U(X,Y).$
	Therefore, $\Sigma$ is totally umbilical if and only if $U$ is a conformal vector field of the degenerate induced metric on $\Sigma$, with conformal factor $-2\mu$ when $B_U = \mu g$. The totally geodesic case is immediate.
\end{proof}

	\subsection{Codimension-Two Spacelike Submanifolds and Factorization}
	For every codimension-two spacelike  submanifold  that factors through a totally umbilical lightlike hypersurface the following result holds.
	
\begin{proposition}\label{ImpliFacil}
	Let $\Sigma$ be a totally umbilical lightlike hypersurface in a spacetime $M$, and let $\psi\colon S \to M$ be a codimension-two spacelike submanifold such that $\psi(S) \subset \Sigma$. Then, there exists an umbilical lightlike normal section $\xi$ of $\psi$. In particular, if $\Sigma$ is a totally geodesic lightlike hypersurface, then $A_{\xi} = 0$ and $\xi$ is a stationary vector field.
\end{proposition}

\begin{proof}
	Consider an arbitrary lightlike vector field $U \in \mathfrak{X}(\Sigma)$ and define $\xi$ as the restriction of $U$ to $S$ along $\psi$, that is, $\xi = U|_{\psi(S)}$. Since $g(X, U) = 0$ for all $X \in \mathfrak{X}(\Sigma)$ and the tangent space of $S$ is contained in the tangent space of $\Sigma$ through $\psi$, it follows that $\xi$ is a lightlike normal section of $S$.
	
	Applying the Weingarten formula~\eqref{fWeingarten}, the definition of the lightlike second fundamental form~\eqref{DefLightlikeSecond}, and the properties of the connection, we obtain
	\begin{equation}\label{EcPropoFactor}
		B_U(X, Y) = g(A_{\xi} X, Y),
	\end{equation}
	for all $X, Y \in \mathfrak{X}(S)$, where $B_U$ is the lightlike second fundamental form of $\Sigma$ associated with $U$, and $A_{\xi}$ is the shape operator corresponding to $\xi$ for the submanifold $S$.
	
	As $\Sigma$ is a totally umbilical lightlike hypersurface, there exists a smooth function $\mu$ on $\Sigma$ such that $B_U = \mu g$. Therefore, from~\eqref{EcPropoFactor}, we have $g(A_{\xi} X, Y) = \mu g(X, Y)$ for all $X, Y \in \mathfrak{X}(S)$. This implies that $A_{\xi} = \rho I$, where $\rho = \mu \circ \psi$ and $I$ is the identity operator. Consequently, $\xi$ is an umbilical lightlike section.
	
	Finally, if $\Sigma$ is a totally geodesic lightlike hypersurface, then $\rho = 0$ and hence $A_\xi = 0$, so that $\xi$ is also a stationary vector field.
\end{proof}

	\begin{example} \normalfont
		This result applies to Brinkmann spacetimes, which are characterized by a globally defined parallel lightlike vector field. These include pp-wave and plane wave spacetimes (\cite{CanovasPalomoRomero} and references therein). Brinkmann spacetimes are foliated by totally geodesic lightlike hypersurfaces, called characteristic hypersurfaces. Thus, any codimension-two spacelike submanifold contained in one of these hypersurfaces has a stationary and umbilical lightlike normal vector field, which can be rescaled to be parallel with respect to the normal connection.
	\end{example}

\begin{proposition} \label{proposition:InterrelationSubmanifoldsHypersurfaces}
	The following interrelations hold between codimension-two spacelike submanifolds and lightlike hypersurfaces in a spacetime $M$:
	\begin{enumerate}
		\item[(a)] Let $\Sigma$ be a lightlike hypersurface in $M$ equipped with an integrable screen distribution $\mathcal{S}$. If $S$ is a leaf of $\mathcal{S}$ and $U$ is a globally defined lightlike vector field on $\Sigma$, then the restriction to $TS$ of the screen shape operator $A^*_U$ (associated with $\mathcal{S}$) coincides with the shape operator $A_{\xi}$ of the spacelike submanifold $S$ in $M$, where $\xi = U|_S$.
		\item[(b)] In the context of the previous item, if the rotation one-form $\tau$ associated with $U$ and ${\mathcal S}$ vanishes, then the induced lightlike normal vector field $\xi = U|_S$ is parallel with respect to the normal connection on $S$ as a submanifold of $M$.
		\item[(c)] If a codimension-two spacelike submanifold $S$ factors through a lightlike hypersurface $\Sigma$, then, regardless of the choice of screen distribution on $\Sigma$, the lightlike expansion scalar of a lightlike vector field $U \in \mathfrak{X}(\Sigma)$ restricted to $S$ is equal to the lightlike expansion scalar of $\xi = U|_S$ on $S$ ($\theta_U|_S = \theta_\xi$), and similarly for the shear scalars ($\sigma_U|_S = \sigma_\xi$).
	\end{enumerate}
\end{proposition}
	\begin{proof}
		\emph{(a)} The result follows directly from~\eqref{RelaLightlik} and~\eqref{EcPropoFactor}. \emph{(b)} This is a consequence of the Weingarten formula~\eqref{fWeingarten} and~\eqref{LightlikeDes1}. 
		\emph{(c)} Regardless of the screen distribution chosen on $\Sigma$, the relation $\theta_U|_S = \theta_\xi$ is established by taking the trace in~\eqref{EcPropoFactor}. Similarly, for shear scalars, it suffices to observe that~$\sigma_U^2|_S = \operatorname{trace}(\mathring{A}^*_U \circ \mathring{A}^*_U)|_S= \operatorname{trace}(\mathring{A}^*_\xi \circ \mathring{A}^*_\xi)=\sigma^2_\xi.$
	\end{proof}

	As an immediate consequence of Propositions~\ref{ImpliFacil}~and~\ref{proposition:InterrelationSubmanifoldsHypersurfaces}, we have the following.
	
\begin{corollary}\label{cor:umbilical_leaves}
	Let $\Sigma$ be a lightlike hypersurface in a spacetime $M$ endowed with an integrable screen distribution $\mathcal{S}$, and let $U$ be a globally defined lightlike vector field on $\Sigma$. Then, $\Sigma$ is totally umbilical if and only if each leaf of $\mathcal{S}$ is umbilical along $U$ as a codimension-two spacelike submanifold of $M$. In particular, $\Sigma$ is totally geodesic if and only if, in addition, the restriction of $U$ to each leaf is stationary.
\end{corollary}
	
	\begin{remark}\normalfont
		For additional general results on the conditions under which a codimension-two spacelike submanifold contained in a lightlike hypersurface is a leaf of the integrable screen distribution, using the rigging technique, see~\cite{gutierrez20}. 
	\end{remark}
	
	\subsection{Conformal Properties and Invariance}
	
The causal character of a hypersurface $\Sigma$, at any point, of a spacetime $(M,g)$ is preserved by pointwise conformal changes of the Lorentzian metric $g$. Furthermore,

\begin{proposition}\label{propo:Conformal_Change}
	Given a lightlike hypersurface $\Sigma$ in $(M,g)$ with lightlike second fundamental form $B$ with respect to a lightlike vector field $U$ on $\Sigma$, if we set $g^* = e^{2u}g$ for $u \in C^{\infty}(M)$, then the corresponding lightlike second fundamental form $B^*$, with respect to $U$, satisfies $B^* = e^{2u}(B - U(u)g)$. In particular, the corresponding lightlike expansion scalar $\theta^*_U$ is given by $\theta^*_U = \theta_U + n U(u)$, where $n = \dim\Sigma-1$. Moreover, if $\Sigma$ is totally umbilical with respect to $g$, then $\Sigma$ is also totally umbilical with respect to $g^*$.
\end{proposition}
\begin{proof}
	If $\widetilde{\nabla}$ and $\widetilde{\nabla}^*$ denote the Levi-Civita connections of $g$ and $g^*$, respectively, then
	\(
	\widetilde{\nabla}^*_X U = \widetilde{\nabla}_X U + X(u)U + U(u)X - g(X,U) \operatorname{grad}_g(u).
	\)
	For $X,Y\in\mathfrak X(\Sigma)$ one has $g(X,U)=g(Y,U)=0$. Hence, using~\eqref{DefLightlikeSecond},
	\[
	\begin{aligned}
	B^*(X,Y)
	&=-g^*(\widetilde{\nabla}^*_XU,Y)\\
	&=e^{2u}\bigl(B(X,Y)-U(u)g(X,Y)\bigr).
	\end{aligned}
	\]
	Thus the screen shape operator changes by a multiple of the identity, namely $A^*_{U,g^*}=A^*_{U,g}-U(u)I$. Since $\theta_U=-\operatorname{trace}(A_U^*)$, the stated transformation law for $\theta_U$ follows. Finally, if $B=\mu g$, then $B^*=(\mu-U(u))g^*$, proving the conformal invariance of total umbilicity.
\end{proof}
	
\begin{remark} \label{CTU} \normalfont
	(a) The lightlike generators of $\Sigma$, which are the integral curves of $U$, transform from lightlike pregeodesics in $(M, g)$ to lightlike pregeodesics in $(M, g^*)$, where $g^* = e^{2u}g$~\cite[Lem.~9.17]{beem}. In particular, if $U$ is a geodesic vector field in $(M, g)$, it becomes a pregeodesic vector field in $(M, g^*)$.
	
	(b) Given a screen distribution $\mathcal{S}$ on $\Sigma$ associated with the metric $g$ of the spacetime $M$, this distribution is also the screen distribution associated with any conformally related metric $g^* = e^{2u}g$. This follows from the fact that the defining properties of a screen distribution are preserved under conformal changes of the metric. Moreover, the integrability of $\mathcal{S}$ depends solely on the closure of the Lie bracket of its vector fields, which is independent of the metric and hence of the conformal transformations. In particular, the integral submanifolds of $\mathcal{S}$ are preserved under conformal changes, so that the foliation defined by $\mathcal{S}$ remains unchanged.
\end{remark}

The following result, which is crucial for our main theorem, describes the behavior of the lightlike expansion scalar and umbilical sections of a codimension-two spacelike submanifold under conformal changes of the ambient Lorentzian metric.

\begin{lemma}\label{lemma:conformal}
	Let $\psi\colon S \to M$ be a codimension-two spacelike submanifold in a spacetime $(M, g)$ that factors through a lightlike hypersurface $\Sigma$, that is, $\psi(S) \subset \Sigma$. Let $U$ be a globally defined lightlike vector field on $\Sigma$, and let $\xi$ denote its restriction to $S$. Then, for any point $p \in S$, there exists a smooth function $u \in C^\infty(M)$ such that, under the conformal change $g^* = e^{2u}g$, the corresponding lightlike expansion scalar $\theta_\xi^*$ is non-zero in a neighborhood of $p$ in $S$.
\end{lemma}

\begin{proof}
	By~\Cref{proposition:InterrelationSubmanifoldsHypersurfaces}(c), the lightlike expansion scalar of $U$ restricted to $S$ coincides with that of $\xi = U|_S$. Under a conformal change $g^* = e^{2u}g$, \Cref{propo:Conformal_Change} gives $\theta_\xi^* = \theta_\xi + n\,\xi(u)$. Thus, for any given point $p \in S$, there exists a smooth function $u$ such that $\theta_\xi^*$ is non-zero in a neighborhood of $p$ in $S$.
\end{proof}

	\subsection{Main Factorization Theorem}
	
The following theorem encapsulates the central geometric construction of this work. This result provides a precise link between the extrinsic geometry of the initial data and the global structure of the associated lightlike hypersurface.

We emphasize that, while a lightlike normal vector field to $S$ always exists locally, the existence of a globally defined lightlike normal vector field is a strong global assumption. The embeddedness of $S$ is also essential, as it guarantees the existence of a well-defined normal exponential map and underpins the global construction in the following theorem.

In our framework, the integrable screen distribution on the constructed lightlike hypersurface is generated by pushing forward the horizontal lifts of vector fields on the initial spacelike submanifold along the lightlike geodesics. This choice makes $S$ itself a leaf and is natural for the geometric analysis developed in this paper.
	
	\begin{theorem}\label{factors}
Let $\psi\colon S\to M$ be a codimension-two spacelike submanifold
embedded in a spacetime $M$, and let
$\xi\in\mathfrak X^\perp(S)$ be a lightlike normal vector field.
Then $\psi$ factors through a lightlike hypersurface $\Sigma$ of $M$
that admits a geodesic lightlike extension
$U\in\mathfrak X(\Sigma)$ of $\xi$ and an integrable screen
distribution $\mathcal S$ for which $S$ is one of the leaves.

If, moreover, $M$ is locally conformally flat and $\xi$ is an
umbilical section of $\psi$, then $\Sigma$ is totally umbilical.
\end{theorem}
	\begin{proof}
	To simplify the notation, we identify $S$ with its image in $M$. Consider the zero section $S_0 = \{ (p, 0) : p \in S \}$ of the normal bundle $T^\perp S$. Since $S$ is embedded, there exists a fiberwise star-shaped open neighborhood $\mathscr V$ of $S_0$ in $T^\perp S$ such that the normal exponential map $\exp^\perp\colon\mathscr V\to\mathscr U:=\exp^\perp(\mathscr V)$ is a diffeomorphism onto a tubular neighborhood $\mathscr U$ of $S$ (see, e.g., \cite[Prop.~7.26]{oneill}). Note that the restriction of $\exp^\perp$ to $S_0$ is simply the diffeomorphism $S_0 \to S$ followed by the embedding $S \hookrightarrow M$.
	
		Let $\xi$ be the lightlike normal section to $S$ given in the statement of the theorem, which we may take, without loss of generality, to be future-directed. Consider the smooth embedding $\Psi \colon S \times \mathbb{R} \to T^\perp S$ defined by $\Psi(p,t) = (p, t\xi_p)$, and set
		\[
			\mathscr O:=\Psi^{-1}(\mathscr V),
			\qquad
			\Phi:=\exp^\perp\circ\Psi|_{\mathscr O}.
		\]
		Then $\mathscr O$ is an open neighborhood of $S\times\{0\}$, and the fiberwise star-shaped property of $\mathscr V$ ensures that
		\[
			I_p:=\{t\in\mathbb R:(p,t)\in\mathscr O\}
		\]
		is an open interval containing $0$ for every $p\in S$. Since $\exp^\perp|_{\mathscr V}$ is a diffeomorphism, $\Phi$ is an embedding and $\Sigma:=\Phi(\mathscr O)\subset\mathscr U$ is an embedded hypersurface of $M$ containing $S$.
		
	We define a vector field $U \in \mathfrak{X}(\Sigma)$ by setting
	\[
	U_{\Phi(p_0,t_0)} = \Phi_{*(p_0,t_0)}\bigl(\partial_t|_{(p_0,t_0)}\bigr)
	\]
	for each $(p_0, t_0) \in \mathscr{O}$, where $\partial_t|_{(p_0,t_0)}$ denotes the tangent vector to the curve $t \mapsto (p_0, t)$ in $S \times \mathbb{R}$ at $t = t_0$. The integral curve of $U$ in $\Sigma$ passing through $\Phi(p_0, t_0) \in \Sigma$ is obtained as the image under the diffeomorphism $\Phi$, which maps $\mathscr{O}$ diffeomorphically onto $\Sigma$, of the integral curve in $\mathscr{O}$ passing through $(p_0, t_0)$. Consequently, it is given by
	\[
	t \mapsto \gamma(t) = \Phi(p_0, t + t_0) = \exp^\perp_{p_0}\bigl((t + t_0)\xi_{p_0}\bigr),
	\]
	which is a geodesic in $M$. Since $\xi$ is lightlike, it is evident that $U$ is also lightlike. Noting that $\Phi(p, 0) = p$ and $U_p = \xi_p$ for all $p \in S \subset \Sigma$, we deduce that $U$ is a geodesic lightlike extension of $\xi$. Moreover, $U$ is future-directed as it coincides with the future-directed $\xi$.
		
		We now establish that $\Sigma$ is a lightlike hypersurface of $M$. Although $\Sigma$ is foliated by lightlike geodesics by construction, this condition alone is not sufficient for $\Sigma$ to be a lightlike hypersurface (see~\Cref{Foliation}). Therefore, we provide a direct argument.
		
		We construct a vector field $X$ on $\Sigma$ that extends a vector field on $S$. Let $X \in \mathfrak{X}(S)$ be an arbitrary vector field on $S$. For any $(p, t) \in \mathscr{O}$, we define $X$ on $\Sigma$ by
		\[
		X_{\Phi(p, t)} = \Phi_{*(p, t)}\bigl(\widehat{X_p}(t)\bigr),
		\]
		where $\widehat{X_p}(t)$ is the horizontal lift of $X_p \in T_p S$ to $(p, t)$. This lift $\widehat{X_p}(t)$ is uniquely determined by requiring that $\pi_{*(p, t)}\bigl(\widehat{X_p}(t)\bigr) = X_p$, where $\pi\colon \mathscr{O} \subseteq S \times \mathbb{R} \to S$ is the natural projection, and that $\widehat{X_p}(t)$ is tangent to $S \times \{t\}$.
		
		Since $\Phi$ is a smooth immersion and the horizontal lift is smooth, this definition produces a smooth vector field $X$ on $\Sigma$. In particular, when $t = 0$, we have $X_p = \Phi_{*(p, 0)}\bigl(\widehat{X_p}(0)\bigr)$ for all $p \in S$. Because $X$ on $\Sigma$ extends $X$ on $S$, we use the same notation $X$ for both vector fields.
		
		The vector field $X$ on $\Sigma$, defined as the pushforward of the horizontal lift, is invariant under the flow $\varphi$ of $U$. To show this, consider $q = \Phi(p,t_0) \in \Sigma$. For any $s$ such that $\varphi_s(q)$ is defined, we have:
		\begin{align*}
			(\varphi_s)_{*q}X_q &= (\varphi_s)_{*\Phi(p,t_0)} \bigl(\Phi_{*(p,t_0)}\bigl(\widehat{X_p}(t_0)\bigr)\bigr) \\
			&= \Phi_{*(p,t_0 + s)}\bigl(\widehat{X_p}(t_0 + s)\bigr) \\
			&= X_{\Phi(p,t_0 + s)} \\
			&= X_{\varphi_s(q)}.
		\end{align*}
		The second equality follows from the relationship between $\varphi$ and $\Phi$. To clarify this, we introduce the translation map $\tau_s\colon S \times \mathbb{R} \to S \times \mathbb{R}$ defined by $\tau_s(p,t) = (p, t + s)$. We then observe that $\varphi_s \circ \Phi = \Phi \circ \tau_s$. Thus:
		\begin{align*}
			(\varphi_s)_{*\Phi(p,t_0)} \bigl(\Phi_{*(p,t_0)}(\widehat{X_p}(t_0))\bigr) &= (\Phi \circ \tau_s)_{*(p,t_0)}\bigl(\widehat{X_p}(t_0)\bigr) \\
			&= \Phi_{*(\tau_s(p,t_0))} \bigl((\tau_s)_{*(p,t_0)}\bigl(\widehat{X_p}(t_0)\bigr)\bigr) \\
			&= \Phi_{*(p,t_0 + s)}\bigl(\widehat{X_p}(t_0 + s)\bigr).
		\end{align*}
		The key step is $(\tau_s)_{*(p,t_0)}\bigl(\widehat{X_p}(t_0)\bigr) = \widehat{X_p}(t_0 + s)$, which holds due to the uniqueness of the horizontal lift. This confirms the invariance of $X$ under $\varphi$.
		
		By the invariance of $X$ under the flow $\varphi$ (see~\cite[Th.~9.42]{lee2}), the vector fields $U$ and $X$ commute, that is, $ \widetilde{\nabla}_U X =\widetilde{\nabla}_X U $. Therefore,
		\begin{align*}
			U\bigl(g(U,X)\bigr) &= g(\widetilde{\nabla}_U U, X) + g(U, \widetilde{\nabla}_U X) \\
			&= g(U, \widetilde{\nabla}_U X) 
			\\
			&= g(U, \widetilde{\nabla}_X U) 
			\\
			&= \frac{1}{2} X\bigl(g(U,U)\bigr) 
			\\
			&= 0.
		\end{align*}
		
		Since $U$ coincides with the lightlike normal vector field $\xi$ along $S$, and $X$ is defined so that $X|_S$ is tangent to $S$, we have $g(U, X) = 0$ at every point of $S$. As $U\bigl(g(U,X)\bigr) = 0$, the function $g(U,X)$ is constant along the integral curves of $U$. Given that the integral curves of $U$ are the lightlike geodesics generating $\Sigma$, it follows that $g(U,X) = 0$ along each lightlike generator. Therefore, $g(U,X) = 0$ for all vector fields $X$ constructed as above. Since these vector fields, together with $U$, span the tangent spaces of $\Sigma$, we conclude that $\Sigma$ is a lightlike hypersurface in $M$.

		The lightlike hypersurface $\Sigma$ admits an integrable screen distribution and is foliated by a family of codimension-two spacelike submanifolds. The distribution $\mathcal{S}$ on $\Sigma$ is defined, for each $(p,t) \in \mathscr{O}$, by
		\[
		\mathcal{S}_{\Phi(p,t)} = \operatorname{span}\bigl\{ \Phi_{*(p,t)}\bigl(\widehat{X_p}(t)\bigr) : X \in {\mathfrak X}(S) \bigr\}.
		\]
		Since $\Phi$ is a diffeomorphism and the horizontal lifts are smooth, the images under $\Phi_{*}$ of the horizontal lifts of smooth local frames on $S$ form smooth local frames for $\mathcal{S}$. This ensures that $\mathcal{S}$ is a smooth distribution on~$\Sigma$.
		
	Note that for each $p \in S$, we have $X_p = \Phi_{*(p,0)}\bigl(\widehat{X_p}(0)\bigr)$, implying that the tangent space $T_pS$ coincides with the subspace $\mathcal{S}_{\Phi(p,0)} = \mathcal{S}_p$ of the distribution $\mathcal S$. It is routine to verify, using the properties of the horizontal lift, pushforward, and Lie brackets, that $\mathcal{S}$ is an involutive distribution complementary to $U$ in $T\Sigma$. By the Frobenius Theorem, $\mathcal{S}$ is an integrable screen distribution on $\Sigma$, and its maximal connected integral manifolds form a foliation of $\Sigma$, with $S$ itself being a leaf of this foliation. This is because $S$ is a connected integral submanifold whose tangent space at each point coincides with the value of the distribution $\mathcal{S}$ at that point. The globally defined leaves obtained at common parameter values are described in~\Cref{remark:submersion}.
		
	Let $\mathfrak{L}(S)$ denote the set of all horizontal lifts of vector fields on $S$ to $\mathscr{O} \subseteq S \times \mathbb{R}$. 
	The set of flow-invariant vector fields obtained as the pushforward of the horizontal lifts of vector fields on $S$ is
	\[
	\mathfrak{I}(\Sigma) = \bigl\{ \Phi_*\bigl(\widehat{X}\bigr) : \widehat{X} \in \mathfrak{L}(S) \bigr\} \subset \Gamma(\mathcal{S}),
	\]
	where $\Gamma(\mathcal{S})$ denotes the set of smooth sections of the distribution $\mathcal{S}$. It is observed that $\mathfrak{I}(\Sigma)$ is a vector subspace of $\Gamma(\mathcal{S})$, but it is not invariant under multiplication by arbitrary smooth functions on~$\Sigma$.
		
	Assume now that $M$ is locally conformally flat and that $\xi$ is an umbilical section. Around any point of $S$, choose a pointwise conformal change $g^*=e^{2u}g$ for which $g^*$ is flat, and identify the resulting neighborhood with an open subset of $\mathbb L^{n+2}$. Umbilicity of the lightlike normal direction is conformally invariant by~\Cref{CIS}. Since the transformed metric is flat, \Cref{CriterionParaUmbi,corollary:ParaConformal} allow us, after a local rescaling, to assume that the transformed lightlike normal field $\xi$ is parallel with respect to the normal connection.

	Write $A_\xi=\rho I$. Since the ambient metric is flat and $\xi$ is parallel, the Codazzi equation gives
	\[
		X(\rho)Y=Y(\rho)X
	\]
	for all $X,Y\in\mathfrak X(S)$; see~\Cref{sec:prelim} for the conventions and~\cite{Dajczerbook} for the standard equation. Since $\dim S\ge2$, it follows that $d\rho=0$. Thus $\rho$ is constant on the connected neighborhood under consideration.

	If $\rho=0$, the submanifold is contained in a degenerate lightlike hyperplane. If $\rho\ne0$, after a constant rescaling and a possible change of sign one may assume that
	\[
		\nabla^\perp\xi=0,
		\qquad
		A_\xi=-I,
	\]
	and the submanifold is contained in a light cone. These classical characterizations are given in~\cite[Theorem~4.3]{IPR}; compare also~\cite[0.9]{Mag}.

	Let $\mathcal H$ denote the degenerate hyperplane or light cone obtained above. In either case, $\xi$ spans $\operatorname{Rad}(\mathcal H)$ along $S$. Hence the lightlike geodesics with initial data $(p,\xi_p)$ are, up to affine reparametrization, precisely the generators of $\mathcal H$ through $p$. By uniqueness of geodesics, after shrinking the domain if necessary, the lightlike hypersurface generated by the normal exponential map coincides with an open subset of $\mathcal H$. In Lorentz--Minkowski spacetime, totally umbilical lightlike hypersurfaces are open subsets of either degenerate hyperplanes or light cones; see~\cite[Corollary~8 and Theorem~13]{AkivisGoldberg}. Thus the generated lightlike hypersurface is totally umbilical for $g^*$.

	A rescaling of $\xi$ changes only the affine parametrization of its lightlike geodesics and does not change the generated hypersurface; see~\Cref{remark:rescaling}. Pointwise conformal changes preserve the unparametrized lightlike generators; see~\Cref{CTU,remark:reparametrize}. Since total umbilicity of a lightlike hypersurface is conformally invariant by~\Cref{propo:Conformal_Change}, $\Sigma$ is totally umbilical on a neighborhood of $S$ for the original metric.

	It remains to propagate this local conclusion along the generators. Fix $p\in S$, write $\gamma_p(t)=\Phi(p,t)$, and let $J_p$ be the set of parameters for which $\Sigma$ is totally umbilical on a neighborhood of $\gamma_p(t)$. The set $J_p$ is non-empty and open. If $t_0$ belongs to its closure, choose a conformally flat neighborhood $\mathscr W$ of $\gamma_p(t_0)$ and $t\in J_p$ sufficiently close to $t_0$. A local leaf $L$ of $\mathcal S$ through $\gamma_p(t)$ is umbilical along $U$ by~\Cref{cor:umbilical_leaves}. The lightlike normal exponential map of $(L,U|_L)$ has the same generators as $\Sigma$ and therefore parametrizes $\Sigma\cap\mathscr W$, after shrinking $\mathscr W$ if necessary. The preceding local argument, applied to $(L,U|_L)$, proves total umbilicity near $\gamma_p(t_0)$. Hence $J_p$ is closed. Since the parameter interval of $\gamma_p$ is connected, $J_p$ is the whole interval. As $p$ was arbitrary, $\Sigma$ is totally umbilical.
	\end{proof}

It is natural to ask whether a totally umbilical lightlike hypersurface containing a stationary umbilical leaf must be totally geodesic. The following example shows that this need not be the case.

\begin{example}\label{counterexample_GRW}
	\normalfont
	Consider a generalized Robertson-Walker spacetime~\cite{AliasRomeroSanchez1995} \( M = (0,\infty) \times_f N \), where \( f \in C^\infty((0,\infty)) \) is a positive function and \( N = (-1,1) \times \mathbb{R}^2 \) is a Riemannian manifold equipped with the metric \( g_N = ds^2 + \delta(s, x, y)^2 (dx^2 + dy^2) \), where \( \delta \) is a positive smooth function on \( N \). The Lorentzian metric on \( M \) is given by \( g = -dt^2 + f(t)^2 g_N \).
	
	For a fixed \( t_0 > 0 \), consider an open interval \( I \subset (0, \infty) \) around \( t_0 \), and let \( \Sigma \) denote the lightlike hypersurface in \( M \) defined by
	\[
	\Sigma = \left\{ (t, s, x, y) \in I \times (-1,1) \times \mathbb{R}^2 : s = \int_{t_0}^t \frac{1}{f(r)}\, dr \right\}.
	\]
	By continuity, the interval \( I \) can be chosen so that, for all \( t \in I \), the value \( s = \int_{t_0}^t \frac{1}{f(r)}\, dr \) remains in \( (-1,1) \). This hypersurface contains the spacelike embedded surface
	\(
	S := \{ (t_0, 0, x, y) : (x, y) \in \mathbb{R}^2 \},
	\)
	which can be identified with \( \mathbb{R}^2 \) at \( t = t_0 \) and \( s = 0 \).
	
	The existence of \( \Sigma \) as a totally umbilical lightlike hypersurface of \( M \), as well as the expression for the lightlike expansion scalar \( \theta = \theta_U \) of \( \Sigma \),
	\[
		\theta = -\frac{2}{f(t)^2} \left( f'(t) + \frac{\delta_s(s, x, y)}{\delta(s, x, y)} \right),
	\]
	follow directly from~\cite[Th.~4.2]{gutierrez15}, which describes this construction in terms of the twisted product decomposition of 
	\( N \). Here, \( \delta_s \) denotes the partial derivative of \( \delta \) with respect to \( s \). Moreover, \( U \in \mathfrak{X}(\Sigma) \) is the unique lightlike vector field satisfying \( g(\mathbf{Z}, U) = 1 \), where \( \mathbf{Z} = f\partial_t \) is a globally defined timelike vector field on \( M \).
	
	Now, define \( \delta(s, x, y) \) on \( N \) by
	\[
	\delta(s, x, y) := (1-s) h(x, y),
	\]
	where \( h \) is a smooth positive function on \( \mathbb{R}^2 \). Clearly, \( \delta \) is smooth and positive, since \( s \in (-1,1) \) ensures \( 1-s > 0 \) and \( h(x, y) > 0 \). Thus, the Riemannian metric \( g_N \) is well-defined.
	
	A straightforward computation shows that
	\[
		\theta = -\frac{2}{f(t)^2} \left( f'(t) - \frac{1}{1-s} \right).
	\]
	Specifically, if we choose \( f(t) = t \), then \( s = \ln(t/t_0) \). For \( s \) to remain in \( (-1,1) \), the interval \( I \) must be chosen such that \( t \in (t_0/e, t_0 e) \). With this choice, we have
	\[
		\theta = \frac{2}{t^2} \left( \frac{s}{1-s} \right),
	\]
	which vanishes only at \( s = 0 \), corresponding to $t=t_0$.
	
	Consequently, taking into account \Cref{proposition:InterrelationSubmanifoldsHypersurfaces}(c), \( S \) is a codimension-two spacelike embedded surface of \( M \) contained in \( \Sigma \), with \( \xi = U|_S \) providing a stationary and umbilical lightlike normal vector field along \( S \). However, while \( \Sigma \) is totally umbilical, it is not totally geodesic.
	
	It is important to note that, locally, the lightlike geodesics emanating from \( S \) in the direction of \( \xi \) remain within \( \Sigma \) (see Remark~\ref{remark:TU_Hyper}). Therefore, in a neighborhood of each point of \( S \), \( \Sigma \) coincides with the lightlike hypersurface constructed via the normal exponential map along lightlike geodesics emanating from \( S \) in the direction of \( \xi \), as described in~\Cref{factors}.
\end{example}

	\subsection{Remarks and Discussion}

	In what follows, we adopt the notation and constructions from the proof of~\Cref{factors}. Specifically, let $S$ be a codimension-two spacelike embedded submanifold in a spacetime $(M, g)$, and let $\xi$ be a lightlike normal section along $S$. We consider the smooth map $\Phi \colon \mathscr{O} \subseteq S \times \mathbb{R} \to M$ defined by the normal exponential map, and set $\Sigma = \Phi(\mathscr{O})$ as the associated lightlike hypersurface containing $S$. We also denote by $U$ the lightlike geodesic extension of $\xi$ along $\Sigma$, and by $\mathcal{S}$ the integrable screen distribution on $\Sigma$ given by the pushforward of horizontal lifts of vector fields on $S$ via $\Phi$. All subsequent statements refer to these objects unless explicitly stated otherwise.

\begin{remark}\label{remark:rescaling}
	\textbf{Invariance under rescaling of the lightlike normal section.} \normalfont
	
	Let $\overline{\xi} = f\xi$ be a rescaling of the lightlike normal section $\xi$ by a non-vanishing smooth function $f$ on $S$. If $\mathscr O$ is the domain used for $\xi$, set
	\[
		\overline{\mathscr O}:=\{(p,s):(p,f(p)s)\in\mathscr O\},
		\qquad
		\overline\Phi(p,s):=\Phi(p,f(p)s).
	\]
	Then $\overline\Phi$ is precisely the normal exponential construction for $\overline\xi$ on $\overline{\mathscr O}$ and has the same image $\Sigma_\xi=\Phi(\mathscr O)$. Thus the two constructions have the same unparametrized lightlike generators and differ only in their affine parametrizations.
	
	The geodesic extension of $\overline{\xi}$ to $\Sigma_\xi$ is given by $\overline{U} = \widetilde{f} U$, where $\widetilde{f}$ denotes the natural extension of $f$ to $\Sigma_\xi$ along the generators, defined by $\widetilde{f}(q) := \bigl(f \circ \pi \circ \Phi^{-1}\bigr)(q)$ for each $q \in \Sigma_\xi$, with $\pi\colon\mathscr{O} \subseteq S \times \mathbb{R} \to S$ the canonical projection.
	
	In particular, geometric properties of $\Sigma_\xi$ that depend only on the direction of $\xi$, such as total umbilicity, are invariant under rescalings of $\xi$ by non-vanishing smooth functions. Therefore, any conformal or geometric property of $\Sigma_\xi$ determined by the direction of $\xi$ remains unchanged under such rescalings.
\end{remark}
	\begin{remark}
		\textbf{Transversal intersection.} \normalfont
		
	Given lightlike normal sections $\xi$ and $\eta$ to $S$ with $g(\xi, \eta) = -1$, \Cref{factors} constructs embedded lightlike hypersurfaces $\Sigma_\xi$ and $\Sigma_\eta$ containing $S$, generated respectively by the normal exponential map along $\xi$ and $\eta$. The condition $g(\xi, \eta) = -1$ ensures that they intersect transversally along $S$. Consequently, for every $p\in S$, after restricting both hypersurfaces to a sufficiently small neighborhood of $p$, their intersection is precisely $S$ in that neighborhood.
		
		This transversal intersection is fundamental, as it enables the study of the geometry of $S$ through the properties of the intersecting lightlike hypersurfaces. If $S$ can be realized locally as the intersection of two totally umbilical lightlike hypersurfaces, then $S$ is totally umbilical by~\Cref{ImpliFacil}. Conversely, when the ambient spacetime is locally conformally flat, \Cref{factors} shows that every totally umbilical codimension-two spacelike submanifold can be realized locally as such an intersection.

	\end{remark}
	
	\begin{remark}\label{remark:cones}
		\textbf{Factorization through light cones.} \normalfont 
		
		In complete spacetimes of constant curvature and dimension greater than three, the classification of totally umbilical lightlike hypersurfaces is particularly rigid. In fact, such a hypersurface must be either totally geodesic or contained in a light cone~\cite{AkivisGoldberg}. Consequently,~\Cref{factors} shows that if a codimension-two spacelike embedded submanifold $S$ in such a spacetime admits an umbilical lightlike normal section $\xi$ with non-vanishing lightlike expansion scalar $\theta_\xi$, then $S$ necessarily factors through a light cone in $M$.

This factorization result has notable implications. For instance, if $S$ is compact and the ambient spacetime is Lorentz-Minkowski, Proposition~5.1 in~\cite{PalmasPalomoRomero} ensures that $S$ must be topologically an $n$-sphere. Similarly, in the closed Friedmann cosmological model, light cones are the only totally umbilical lightlike hypersurfaces~\cite{gutierrez15}. Thus,~\Cref{factors} implies that, in this cosmological context, any codimension-two spacelike embedded submanifold with an umbilical lightlike normal section must also factor through a light cone.
		
	\end{remark}
		\begin{remark}\label{remark:reparametrize}
		\textbf{Invariance of equivalence classes of flow-invariant vector fields under pointwise conformal changes of the ambient metric.} \normalfont
		
	(a) Let $g^* = e^{2u}g$ denote a pointwise conformal rescaling of the ambient Lorentzian metric for some $u \in C^\infty(M)$. 
	Let $\overline{U}$ be the global rescaling of $U$ on $\Sigma$ such that $\overline{U}$ is geodesic with respect to the Levi-Civita connection associated to $g^*$ (see Remark~\ref{CTU}(a)). 
	In this setting, $\overline{U}$ must be of the form
	\[
	\overline{U} = f U, \qquad f = C\,e^{-2u},
	\]
	where $C$ is a non-vanishing smooth function that is constant along each lightlike generator of $\Sigma$. 
	This explicit form of $f$ is obtained by requiring that $\overline{U}$ be geodesic for $g^*$, using the standard relation between the Levi-Civita connections of $g$ and $g^*$ under conformal changes; see the proof of \Cref{propo:Conformal_Change}.
	
	As in the proof of~\Cref{factors}, all constructions can be carried out analogously by replacing $g$ with $g^*$ and $U$ with $\overline{U}$. In particular, the integrable screen distribution constructed from the flow-invariant vector fields associated with $\overline{U}$ and $g^*$ is, at each $q \in \Sigma$, given by 
	\[
	\overline{\mathcal{S}}_q = \operatorname{span}\bigl\{ \overline{\Phi}_{*(p, s)}\bigl(\widehat{X_p}(s)\bigr) : X \in \mathfrak{X}(S),\ \overline{\Phi}(p, s) = q \bigr\},
	\]
	where $\overline{\Phi}\colon \overline{\mathscr{O}} \subseteq S \times \mathbb{R} \to \Sigma$ parametrizes the flow lines of $\overline{U}$, and the horizontal lift and pushforward are defined as in the original proof, but with respect to $g^*$.
	Thus, regardless of any pointwise conformal change of the ambient Lorentzian metric, both $\mathcal{S}$ and $\overline{\mathcal{S}}$ define integrable screen distributions on $\Sigma$, as discussed in Remark~\ref{CTU}(b).
	
	Now, let $\Phi\colon \mathscr{O} \subseteq S \times \mathbb{R} \to \Sigma$ be the diffeomorphism constructed in the proof of~\Cref{factors} corresponding to $U$ and $g$. 
	For each $p \in S$, the curve $\gamma(t) = \Phi(p, t)$ is the lightlike generator with affine parameter $t$ and initial conditions $\gamma(0) = p$ and $\gamma'(0) = U_p$. 
	Similarly, $\overline{\gamma}(s) = \overline{\Phi}(p, s)$ is the lightlike generator for $\overline{U}$, with $\overline{\gamma}(0) = p$ and $\overline{\gamma}'(0) = \overline{U}_p$. 
	
	These two curves trace the same geometric path in $\Sigma$, differing only by a reparametrization. 
	For each $p \in S$, there exists a strictly monotonic smooth function $h_p\colon \overline{I}_p \to I_p$ such that $\overline{\gamma}(s) = \gamma(h_p(s))$ and $h_p(0)=0$, where
	\[
	I_p = \{ t \in \mathbb{R} : (p, t) \in \mathscr{O} \}, \qquad \overline{I}_p = \{ s \in \mathbb{R} : (p, s) \in \overline{\mathscr{O}} \}
	\]
	are the open intervals of definition for the respective parametrizations along the lightlike generator through $p$. Define $h\colon \overline{\mathscr{O}} \to \mathbb{R}$ by $h(p, s) := h_p(s)$, so that the domain of $h$ is $\overline{\mathscr{O}} \subseteq S \times \mathbb{R}$. 
	For each fixed $s$, let $h_s(p) := h(p, s)$, defined on $\overline{\mathscr{O}}_s := \{ p \in S : (p, s) \in \overline{\mathscr{O}} \}$. Note that $\overline{\mathscr{O}}_s$ is an open subset of $S$.
	
To relate the screen distributions $\mathcal{S}$ and $\overline{\mathcal{S}}$ of $\Sigma$, let $X \in \mathfrak{X}(S)$, and consider the generators $X_q$ and $\overline{X}_q$ of the distributions at $q = \Phi(p, t) = \overline{\Phi}(p, s) \in \Sigma$, constructed as $X_q = \Phi_{*(p, t)}\bigl(\widehat{X_p}(t)\bigr)$ and $\overline{X}_q = \overline{\Phi}_{*(p, s)}\bigl(\widehat{X_p}(s)\bigr)$. Fix $p \in S$ and, for each fixed $s$, consider a smooth curve $\alpha$ in $S$ with $\alpha(0) = p$ and $\alpha'(0) = X_p$, such that $\alpha(\varepsilon) \in \overline{\mathscr{O}}_s$ for all $\varepsilon$ in a sufficiently small interval around $0$. Using the Chain Rule and the relation $\overline{\Phi}\bigl(\alpha(\varepsilon), s\bigr) = \Phi\bigl(\alpha(\varepsilon), h_s(\alpha(\varepsilon))\bigr)$, we compute
\begin{align*}
	\overline{X}_q &= \frac{d}{d\varepsilon}\bigg|_{\varepsilon=0} \overline{\Phi}\bigl(\alpha(\varepsilon), s\bigr) \\
	&= \frac{d}{d\varepsilon}\bigg|_{\varepsilon=0} \Phi\bigl(\alpha(\varepsilon), h_s(\alpha(\varepsilon))\bigr) \\
	&= \Phi_{*(p,t)} \bigl(X_p, k\bigr),
\end{align*}
where $t = h_p(s)$ and $k = X_p(h_s)$. Here, under the standard identification $T_{(p, t)}(S \times \mathbb{R}) \simeq T_pS \oplus \mathbb{R}$, the vector $(X_p, 0)$ corresponds to the horizontal lift $\widehat{X_p}(t)$ and $(0, 1)$ to the coordinate vector field $\partial_t$ at $(p, t)$. For each $(p, s)\in \overline{\mathscr{O}}$, $k = X_p(h_s)$ is the directional derivative of the function $h_s$ at the point $p$ in the direction $X_p$. In particular, for $s = 0$, $h_0(p) = 0$ for all $p$, so $k = 0$ along $S$.

Now, by linearity of the pushforward and using the previous expression, we can write
\begin{align}
	\overline{X}_q &= \Phi_{*(p, t)}\left((X_p, 0) + (0, k)\right) \nonumber \\
	&= \Phi_{*(p, t)}\bigl(\widehat{X_p}(t)\bigr) + k\, \Phi_{*(p, t)}\bigl(\partial_t\bigr) \nonumber \\
	&= X_q + k\, U_q, \label{eq:Xbar_X_relation}
\end{align}
where $k = X_p(h_s)$. 
	
Therefore, although $\mathcal{S}_q$ and $\overline{\mathcal{S}}_q$ are, in general, distinct subspaces of $T_q\Sigma$, their corresponding flow-invariant generators $X_q$ and $\overline{X}_q$ differ only by lightlike multiples in $\operatorname{Rad}_q(\Sigma)$. This means that the equivalence classes $[X_q]$ and $[\overline{X}_q]$ in the quotient space $T_q\Sigma / \operatorname{Rad}_q(\Sigma)$ coincide for all $q \in \Sigma$. Thus, for each \( X \in \mathfrak{X}(S) \), the equivalence class of the associated flow-invariant vector field in the quotient bundle \( T\Sigma / \operatorname{Rad}(\Sigma) \) remains unchanged under any pointwise conformal changes of the ambient Lorentzian metric.

(b) We note that relation~\eqref{eq:Xbar_X_relation} also holds if, from the outset, one considers an arbitrary non-vanishing smooth rescaling $\overline{U} = fU$ on $\Sigma$, not necessarily arising from a pointwise conformal change of the ambient Lorentzian metric. In this situation, the construction of the integrable screen distribution $\overline{\mathcal{S}}$ associated to $\overline{U}$ follows by an analogous argument to the proof of~\Cref{factors}, applied to the rescaled (pregeodesic) vector field. The explicit formula~\eqref{eq:Xbar_X_relation} follows by the same reasoning as above, and the corresponding vector fields $\overline{X}$ are invariant under the flow of $\overline{U}$, as in the original construction.

In particular, in this context, formula~\eqref{eq:Xbar_X_relation} shows that, although the integrable screen distribution constructed from flow-invariant vector fields (and thus the smooth structure of its integral spacelike submanifolds, except for $S$ itself) depends on the specific choice of lightlike vector field in the radical distribution, its equivalence class in the quotient by the radical remains invariant under such rescalings. This invariance is closely related to the canonical viewpoint of Kupeli, as discussed in Remark~\ref{remark:canonical_screen_Kupeli}.
	\end{remark}

	\begin{remark}\label{remark:submersion}
		\textbf{Associated submersion structure.} \normalfont

		The smooth map
		\[
			F\colon\Sigma\longrightarrow S,
			\qquad
			F=\pi\circ\Phi^{-1},
		\]
		is, by construction, a surjective smooth submersion. Its fibers
		$F^{-1}(p)$, $p\in S$, are the lightlike geodesic generators of
		$\Sigma$ issuing from $p$. The screen distribution $\mathcal S$ is
		orthogonal to these fibers and, for every
		$X\in\mathfrak I(\Sigma)$ and $q=\Phi(p,t)\in\Sigma$,
		\[
			\mathrm dF_q(X_q)=X_p.
		\]
		Thus the flow-invariant vector fields are precisely the horizontal
		lifts, with respect to $F$, of vector fields on $S$~\cite{oneill}.

		For each $t\in\mathbb R$, set
		\[
			\mathscr O_t:=\{p\in S:(p,t)\in\mathscr O\}.
		\]
		Since the screen distribution constructed in~\Cref{factors} is the
		pushforward by $\Phi$ of the distribution tangent to the slices
		$S\times\{t\}$, every leaf $L$ of $\mathcal S$ is the image under
		$\Phi$ of a connected component of
		$\mathscr O_t\times\{t\}$ for a unique $t$. Consequently,
		\[
			F|_L\colon L\longrightarrow F(L)
		\]
		is a diffeomorphism, and $F(L)$ is a connected component of
		$\mathscr O_t$. In general, $F(L)$ need not be all of $S$.
	\end{remark}

	\begin{remark}\label{remark:common-domain}
		\textbf{Common parameter domain and uniformity.} \normalfont

		For each $p\in S$, let
		\[
			I_p:=\{t\in\mathbb R:(p,t)\in\mathscr O\},
			\qquad
			D:=\bigcap_{p\in S}I_p.
		\]
		Each $I_p$ is an open interval containing $0$, because $\mathscr O$
		is open and fiberwise star-shaped. Hence $D$ is an interval containing
		$0$, although it need not be open. For every $t\in D$, the map
		$\Phi_t(p):=\Phi(p,t)$ is an embedding of $S$, its image
		$S_t:=\Phi_t(S)$ is a leaf of $\mathcal S$, and
		$F|_{S_t}\colon S_t\to S$ is a diffeomorphism with inverse $\Phi_t$.
		We refer to the $S_t$, $t\in D$, as the globally defined
		\emph{canonical leaves}.

		The underlying distinction is the familiar one between a tubular
		neighborhood and a uniform tubular neighborhood. O'Neill's
		normal-neighborhood theorem~\cite[Prop.~7.26]{oneill}, used in the
		proof of~\Cref{factors}, provides an open neighborhood of the zero
		section on which the normal exponential map is a diffeomorphism, but
		it does not impose a common fiber size over a non-compact base. In
		the Riemannian setting, Lee makes this distinction explicit:
		every embedded submanifold admits a tubular neighborhood, whereas a
		compact embedded submanifold admits a uniform tubular neighborhood
		\cite[Th.~5.25]{LeeRiem}.
	\end{remark}

	\begin{lemma}\label{lemma:compact-common-domain}
		If $S$ is compact, then $D$ is an open interval containing $0$.
	\end{lemma}

	\begin{proof}
		Let $t_0\in D$. Then $S\times\{t_0\}\subset\mathscr O$. Since $S$ is
		compact and $\mathscr O$ is open in $S\times\mathbb R$, the tube
		lemma~\cite[Lem.~26.8]{Munkres}, applied after interchanging the two
		factors, gives an open interval $J$ containing $t_0$ such that
		$S\times J\subset\mathscr O$. Hence $J\subset D$, so $D$ is open.
		Since $D$ is an intersection of intervals containing $0$, it is itself
		an interval containing $0$.
	\end{proof}

	\begin{remark}\normalfont
		Even when $S$ is compact, the family of canonical leaves need not
		exhaust $\Sigma$; equivalently, $\Phi(S\times D)$ need not equal
		$\Sigma$.
	\end{remark}

	\begin{example}\label{example:nonuniform-leaf-domain}
		\normalfont
		The failure of global surjectivity already occurs in Lorentz--Minkowski
spacetime. Let $n\ge2$ and consider $\mathbb L^{n+2}$ with coordinates
$(u,x_1,\ldots,x_n,v)$ and metric
\(g=-du^2+dx_1^2+\cdots+dx_n^2+dv^2\).
For \(x=(x_1,\ldots,x_n)\in\mathbb R^n\), let
\[
 r(x)=\frac{1}{1+(\|x\|^2-1)^2},
\]
where $\|\cdot\|$ denotes the Euclidean norm on $\mathbb R^n$, and consider
the open spacetime
\(
M=\{(u,x,v)\in\mathbb L^{n+2}:|u|<r(x),\ |v|<r(x)\}.
\)
Let \(S=\{(0,x,0):x\in\mathbb R^n\}\) and take the constant lightlike
normal vector field \(\xi=\partial_u+\partial_v\) along $S$. Then $S$ is
an embedded spacelike submanifold of $M$, and the normal exponential
construction in the direction of $\xi$ is
\[
 \mathscr O
 =
 \{(x,t)\in\mathbb R^n\times\mathbb R:|t|<r(x)\},
 \qquad
 \Phi(x,t)=(t,x,t).
\]
Thus $\Phi$ is an embedding and distinct lightlike generators do not
intersect.

For $t=1/2$,
\(
\mathscr O_{1/2}
=\{x\in\mathbb R^n:r(x)>1/2\}
=\{x\in\mathbb R^n:0<\|x\|^2<2\},
\)
which is a connected punctured ball. Hence, for the leaf
\(L=\Phi(\mathscr O_{1/2}\times\{1/2\})\),
\[
 F(L)=\mathscr O_{1/2}
 =\{x\in\mathbb R^n:0<\|x\|^2<2\}
 \subsetneq S\simeq\mathbb R^n.
\]
Therefore $F|_L\colon L\to F(L)$ is a diffeomorphism, but $F|_L$ is not
surjective onto $S$. In particular, $L$ is not diffeomorphic to $S$, since
\(H_{n-1}(L;\mathbb Z)\cong\mathbb Z\), whereas
\(H_{n-1}(\mathbb R^n;\mathbb Z)=0\).
	\end{example}

\begin{remark}\label{remark:local_version}
		\textbf{Local applicability of the Factorization Theorem.} \normalfont
		
		\Cref{factors} (the Factorization Theorem) provides a local characterization for any totally umbilical lightlike hypersurface $\Sigma$ in a spacetime $M$. Indeed, lightlike geodesics of $M$ with initial vectors tangent to $\Sigma$ remain geodesic in $\Sigma$ locally (see~\Cref{remark:TU_Hyper}). The local nature of the result is ensured by~\Cref{ImpliFacil}, which shows that any totally umbilical lightlike hypersurface locally contains a codimension-two spacelike submanifold $S$ that admits an umbilical lightlike normal section. Consequently, all results discussed above, as well as those presented in the next section, admit local adaptations for any totally umbilical lightlike hypersurface in a spacetime.
		
		This highlights the strength and flexibility of the factorization framework, making it a powerful tool for both local and global studies of totally umbilical lightlike hypersurfaces in Lorentzian geometry.
	\end{remark}
	
	\medskip
	
Although the containment of a codimension-two spacelike embedded submanifold in a lightlike hypersurface can be established by classical methods, such as those employing Fermi coordinates and the normal exponential map~\cite[p.~60]{penrose1972}, our approach offers several distinct advantages. By explicitly constructing globally defined invariant vector fields and a screen distribution, we obtain direct access to important geometric structures associated with the hypersurface.

The normal exponential construction provides a lightlike hypersurface with a geodesic lightlike extension and an integrable screen distribution. By~\Cref{cor:umbilical_leaves}, the generated hypersurface is totally umbilical precisely when every leaf of this screen distribution is umbilical along the restriction of the lightlike generator. In particular, the umbilicity of the initial lightlike normal section propagates whenever the ambient spacetime is locally conformally flat.

	\section{Consequences of the Factorization Theorem}\label{sec:conse_Umbil_Fac_Th}

	The Factorization Theorem establishes a deep link between embedded codimension-two spacelike submanifolds and lightlike hypersurfaces. Throughout this section, whenever an umbilical conclusion is used, we assume explicitly that the generated lightlike hypersurface $\Sigma$ is totally umbilical. By~\Cref{cor:umbilical_leaves}, this is equivalent to requiring that every leaf of the associated screen distribution be umbilical along $U$. When the ambient spacetime is locally conformally flat, \Cref{factors} shows that this hypothesis follows from the umbilicity of the initial lightlike normal section $\xi$.
	
	\subsection{Conformal Relationship}
	
	If $\Sigma$ is totally umbilical,~\Cref{lemmaULH} shows that any lightlike vector field $U$ on $\Sigma$ is conformal with respect to the induced degenerate metric. Although general aspects of conformal structures on lightlike hypersurfaces have been explored in works such as~\cite{Duggal_Sahin,DuggalSharma}, the submersion $F$ and the normal exponential parametrization allow us to compare explicitly the induced Riemannian metrics on those leaves that project globally onto the initial leaf $S$:
	
	\begin{theorem}\label{Cor:ConformalRelation}
		Let $S$ be a spacelike codimension-two embedded submanifold of a
		spacetime $M$ with a lightlike normal vector field $\xi$. Let
		$\Sigma$ be the lightlike hypersurface emanating from $S$ via the
		map $\Phi$ of~\Cref{factors}, and assume that $\Sigma$ is totally
		umbilical. Let $F\colon\Sigma\to S$ be the submersion introduced
		in~\Cref{remark:submersion}. Then:
		\begin{enumerate}
			\item[(a)] For any $X,Y\in\mathfrak X(\Sigma)$ and each
			$q=\Phi(p,t)\in\Sigma$, the induced degenerate metric satisfies
			\begin{equation}\label{Eq:ConformalProperUmbilical_0}
				g_q(X_q,Y_q)
				=\Omega(q)\,
				g_p\bigl(\mathrm dF_q(X_q),\mathrm dF_q(Y_q)\bigr),
			\end{equation}
			where
			\begin{equation}\label{Eq:ConformalRemarkProperUmbilical}
				\Omega\bigl(\Phi(p,t)\bigr)
				=\exp\left(-2\!\int_0^t
				\mu\bigl(\Phi(p,s)\bigr)\,ds\right),
			\end{equation}
			and $\mu$ is determined by $B_U=\mu g$ for the geodesic
			extension $U$ of $\xi$.

			\item[(b)] Let $L$ and $\overline L$ be two leaves of
			$\mathcal S$ such that
			$F|_L\colon L\to S$ and
			$F|_{\overline L}\colon\overline L\to S$ are diffeomorphisms.
			Then
			\[
				G:=(F|_{\overline L})^{-1}\circ F|_L
				\colon L\longrightarrow\overline L
			\]
			is conformal. Consequently, the Riemannian metrics induced on
			$L$ and $\overline L$ are globally conformally related.
			In particular, this applies to any two canonical leaves
			$S_t$ and $S_s$ with $s,t\in D$.
		\end{enumerate}
	\end{theorem}
	\begin{proof}
		\emph{(a)} Let $U=\Phi_*(\partial_t)$ be the geodesic extension of
		$\xi$ to $\Sigma$. Since $\Sigma$ is totally umbilical,
		\Cref{lemmaULH} gives $\mathscr L_Ug=-2\mu g$.

		Fix $p\in S$ and $v,w\in T_pS$, and choose flow-invariant
		extensions $X,Y\in\mathfrak I(\Sigma)$ with $X_p=v$ and $Y_p=w$.
		Along the generator $r\mapsto\Phi(p,r)$ one has
		$[U,X]=[U,Y]=0$, and therefore
		\[
			\frac{d}{dr}g(X,Y)_{\Phi(p,r)}
			=-2\mu\bigl(\Phi(p,r)\bigr)
			g(X,Y)_{\Phi(p,r)}.
		\]
		Integration yields
		\[
			g(X,Y)_{\Phi(p,t)}
			=\Omega\bigl(\Phi(p,t)\bigr)g_p(v,w).
		\]
		Now let $q=\Phi(p,t)$ and take arbitrary
		$X_q,Y_q\in T_q\Sigma$. The restriction of $\mathrm dF_q$ to
		$\mathcal S_q$ is an isomorphism onto $T_pS$, while the remaining
		components are multiples of $U_q$. Since $U_q$ spans the radical
		of $g_q$, the preceding identity applied to the corresponding
		flow-invariant horizontal lifts gives
		\eqref{Eq:ConformalProperUmbilical_0}.

		\emph{(b)} Let $L$ be a leaf for which
		$F|_L\colon L\to S$ is a diffeomorphism. For each $p\in S$, let
		$t_p$ be the unique parameter such that $\Phi(p,t_p)\in L$.
		The function $p\mapsto t_p$ is smooth because
		$(F|_L)^{-1}$ and $\Phi^{-1}$ are smooth. Restricting part~(a)
		to $L$ gives
		\[
			(F|_L)^*g_S
			=e^{2(\phi\circ F|_L)}g_L,
			\qquad
			\phi(p)=\int_0^{t_p}\mu\bigl(\Phi(p,r)\bigr)\,dr.
		\]
		For $\overline L$, define analogously $\overline t_p$ and
		$\overline\phi$. Then
		\[
			(F|_{\overline L})^*g_S
			=e^{2(\overline\phi\circ F|_{\overline L})}
			g_{\overline L}.
		\]
		Therefore
		\begin{equation}\label{eq:h_confromal}
			G^*g_{\overline L}
			=e^{2((\phi-\overline\phi)\circ F|_L)}g_L,
			\qquad
			(\phi-\overline\phi)(p)
			=\int_{\overline t_p}^{t_p}
			\mu\bigl(\Phi(p,r)\bigr)\,dr .
		\end{equation}
		This is precisely the conformal comparison used in the original
		argument. Thus $G$ is conformal. For $L=S_t$ and
		$\overline L=S_s$, the hypotheses follow from
		\Cref{remark:common-domain}.
	\end{proof}

\begin{proposition}  \label{corollary:invariance_omega_en}
		Under the assumptions of~\Cref{Cor:ConformalRelation}, the conformal factor $\Omega$ is invariant under non-vanishing smooth rescalings of the lightlike generator. More precisely, if $\overline U=fU$ and $\overline\Phi(p,s)=\Phi(p,h_p(s))$, then
		\[
			\overline\Omega\bigl(\overline\Phi(p,s)\bigr)
			=\Omega\bigl(\Phi(p,h_p(s))\bigr)
		\]
		whenever both sides are defined.
	\end{proposition}

\begin{proof}
	Let $\overline{U} = f U$, where $f \colon \Sigma \to \mathbb{R}$ is a non-vanishing smooth function (not necessarily positive). As noted in~\Cref{remark:rescaling}, the hypersurface $\Sigma$ remains unchanged under this rescaling, although the parameterizations of the integral curves generally differ. 
	
	The explicit relation between the parameterizations of the integral curves of $U$ and those of $\overline{U}$ is described in~\Cref{remark:reparametrize}(b): for each $p \in S$, there exists a smooth and strictly monotonic function $h_p(s)$ such that $\overline{\gamma}(s) = \gamma(h_p(s))$, where $\gamma(t) = \Phi(p, t)$ is the geodesic integral curve of $U$, and $\overline{\gamma}(s) = \overline{\Phi}(p, s)$ is the (generally pregeodesic) integral curve of $\overline{U}$ starting at $p$ (see \Cref{remark:reparametrize}(a) for details).
	
	Since $\Sigma$ is totally umbilical with respect to $\overline{U}$, we have $B_{\overline{U}} = \overline{\mu}\, g$ for some smooth function $\overline{\mu}$. The relation $B_{\overline{U}} = f B_U$ then gives $\overline{\mu} = f \mu$.
	Although $\overline U$ need not be geodesic, the argument leading to~\eqref{Eq:ConformalRemarkProperUmbilical} uses only the flow equation $\mathscr L_{\overline U}g=-2\overline\mu g$ and therefore applies unchanged to this pregeodesic parametrization.
	
	For each $p \in S$, the reparametrization $t = h_p(s)$ satisfies
	\[
	\frac{dh_p}{ds} = f\bigl(\Phi(p, h_p(s))\bigr),
	\]
	with $h_p(0) = 0$. This follows by differentiating $\overline{\gamma}(s) = \gamma(h_p(s))$.
	
	Therefore, from~\eqref{Eq:ConformalRemarkProperUmbilical}, the conformal factor associated with $\overline{U}$ is
	\[
	\overline{\Omega}\bigl(\overline{\Phi}(p, s)\bigr) = \exp\left(-2 \int_0^s \overline{\mu}\bigl(\overline{\Phi}(p, s')\bigr)\, ds'\right)
	= \exp\left(-2 \int_0^s f\bigl(\overline{\Phi}(p, s')\bigr)\, \mu\bigl(\overline{\Phi}(p, s')\bigr)\, ds'\right).
	\]
	Making the change of variable $t' = h_p(s')$, so that $ds' = dt'/f\bigl(\Phi(p, t')\bigr)$ and $h_p(0) = 0$, $h_p(s) = t$, we obtain
	\[
	\int_0^s f\bigl(\overline{\Phi}(p, s')\bigr)\, \mu\bigl(\overline{\Phi}(p, s')\bigr)\, ds' = \int_0^{h_p(s)} \mu\bigl(\Phi(p, t')\bigr)\, dt'.
	\]
	Therefore, at corresponding points of $\Sigma$,
	\[
	\overline{\Omega}\bigl(\overline{\Phi}(p, s)\bigr) = \Omega\bigl(\Phi(p, h_p(s))\bigr),
	\]
	for any $(p, s)$ in the domain of $\overline{\Phi}$. This proves the assertion.
\end{proof}

\begin{remark}
	\label{remark:Omega_conformal_invariance}
	\normalfont
	
	Under a pointwise conformal change of the ambient metric, $g^* = e^{2u}g$ for some $u \in C^{\infty}(M)$, let $\Omega$ and $\Omega^*$ denote the conformal factors determined by~\eqref{Eq:ConformalProperUmbilical_0} with respect to $g$ and $g^*$, respectively. If $q\in\Sigma$ and $p=F(q)$, then
	\[
		\Omega^*(q)=e^{2(u(q)-u(p))}\Omega(q).
	\]
	Indeed, $g_q^*=e^{2u(q)}g_q$, whereas the metric induced on the initial leaf satisfies $g_p^*=e^{2u(p)}g_p$; the formula follows immediately from~\eqref{Eq:ConformalProperUmbilical_0}. Thus the conformal relationship in~\Cref{Cor:ConformalRelation}(b) is preserved by conformal changes of the ambient metric; in particular, this applies to the canonical leaves $S_t$, $t\in D$. The particular conformal factor, however, transforms according to the preceding law.
\end{remark}
	
In general, the globally defined leaves $S_t$, $t\in D$, need not be isometric. The following corollaries make their conformal relationship precise.
	
\begin{corollary}\label{Cor:isometry_condition}
		Under the assumptions of~\Cref{Cor:ConformalRelation}, let $L$ and
		$\overline L$ be two leaves such that
		$F|_L\colon L\to S$ and
		$F|_{\overline L}\colon\overline L\to S$ are diffeomorphisms.
		With the notation used in the proof of~\Cref{Cor:ConformalRelation},
		the canonical conformal diffeomorphism
		$G=(F|_{\overline L})^{-1}\circ F|_L$ is an isometry if and only if
		\[
			\int_{\overline t_p}^{t_p}
			\mu\bigl(\Phi(p,r)\bigr)\,dr=0
		\]
		for every $p\in S$.
	\end{corollary}
	\begin{proof}
		This is immediate from~\eqref{eq:h_confromal}.
	\end{proof}

	\begin{remark}\normalfont
		If $\Sigma$ is totally geodesic, then $\mu=0$ and
		\Cref{Cor:isometry_condition} shows that any two leaves satisfying
		the global hypothesis of~\Cref{Cor:ConformalRelation}(b) are
		canonically isometric. In particular, all canonical leaves
		$S_t$, $t\in D$, are canonically isometric.

		Moreover, $A_{U|_{S_t}}=0$, so
		\Cref{remark:stationary_umbilical}(b) gives
		$g(\mathbf H_t,\mathbf H_t)=0$, where $\mathbf H_t$ is the mean
		curvature vector of $S_t$. By~\eqref{ExtrinsicScalarH}, the
		extrinsic scalar curvature of every $S_t$ also vanishes.
	\end{remark}

	We finally consider the stronger condition that $\Omega$ itself be
	constant on $\Sigma$.

\begin{corollary}\label{Cor:Omega_cte}
		Under the assumptions of~\Cref{Cor:ConformalRelation}, if the
		conformal factor $\Omega$ is constant on $\Sigma$, then
		$\Omega\equiv1$ and $\Sigma$ is totally geodesic. Consequently,
		every leaf $L$ for which $F|_L\colon L\to S$ is a diffeomorphism
		is canonically isometric to $S$; in particular, this holds for every
		canonical leaf $S_t$, $t\in D$.
	\end{corollary}
	\begin{proof}
		Since $\Omega|_S=1$, constancy gives $\Omega\equiv1$.
		Differentiating~\eqref{Eq:ConformalRemarkProperUmbilical} along the
		generators yields $\mu=0$ on $\Sigma$, and hence $\Sigma$ is
		totally geodesic. The conclusion follows from
		\Cref{Cor:isometry_condition}.
	\end{proof}

\subsection{Volume Evolution of Compact Leaves}

Let \(M\) be an oriented spacetime. Then, all codimension-two spacelike submanifolds with a lightlike normal are orientable (see~\Cref{notaOrientable}(a)). In particular, they possess a well-defined volume.

In the framework of~\Cref{factors}, let \(S\) be a compact spacelike submanifold with a lightlike normal vector field \(\xi\), and assume that the lightlike hypersurface \(\Sigma\) generated by the normal exponential map is totally umbilical. By~\Cref{lemma:compact-common-domain}, \(D\) is an open interval containing \(0\), and \(\{S_t\}_{t\in D}\) is a family of compact spacelike leaves of the screen distribution, with \(S_0=S\). Each \(S_t\) has a well-defined volume. Recall that \(\Phi(S\times D)\) need not exhaust \(\Sigma\).

We now study how the volume of the spacelike submanifolds \(S_t\) evolves as \(t\) varies, focusing on the relationship between the volume element on \(S\) and the volume element on \(S_t\) induced by the ambient Lorentzian metric.

\begin{remark}\normalfont
	For any $t_0\in D$, the same argument applies with $S_{t_0}$ and its induced volume density as the initial data, and with the integral along each generator taken from $t_0$ to $t$.
\end{remark}

\begin{theorem}\label{Th:Vol1}
	Let \(S\) be a compact spacelike codimension-two embedded submanifold of an orientable spacetime \(M\), admitting a lightlike normal vector field \(\xi\). Let \(\Sigma\) be the lightlike hypersurface emanating from \(S\) through the map \(\Phi\) defined in~\Cref{factors}, and assume that \(\Sigma\) is totally umbilical. Then, for each \(t \in D\), the volume of a leaf \(S_t\) of the screen distribution \(\mathcal{S}\), induced by the metric of \(M\), is given by
	\begin{equation}\label{Vol1}
		\operatorname{Vol}(S_t) = \int_S \exp\left( \int_0^t \theta\bigl(\Phi(p, s)\bigr) \, ds \right) {\rm d}V,
	\end{equation}
	where \(\theta\) denotes the lightlike expansion scalar associated with the lightlike geodesic extension \(U\) of \(\xi\) to \(\Sigma\), and \({\rm d}V\) is the induced volume element on \(S\).
\end{theorem}
\begin{proof}
	Fix the lightlike geodesic extension \(U\) of \(\xi\) and its affine parameter \(t\). Denote the induced volume elements on \(S\) and \(S_t\) by \({\rm d}V\) and \({\rm d}V_t\), respectively. By~\Cref{Cor:ConformalRelation}, the diffeomorphism \(\Phi_t\colon S\to S_t\) satisfies
	
	\[
	\Phi_t^* {\rm d}V_t = (\Omega \circ \Phi_t)^{n/2} \, {\rm d}V,
	\]
	where \(n = \dim S\).
	
	Hence, for each \(t \in D\), the volume of \(S_t\) is
	\[
	\operatorname{Vol}(S_t) = \int_S (\Omega \circ \Phi_t)^{n/2} \, {\rm d}V.
	\]
	Substituting \((\Omega \circ \Phi_t)(p) = \exp\left(-2 \int_0^t \mu(\Phi(p, s)) \, ds\right)\) from~\Cref{Cor:ConformalRelation}, and using \(\theta = -n \mu\) for the totally umbilical lightlike hypersurface \(\Sigma\), we obtain the volume evolution formula~\eqref{Vol1}.
\end{proof}

\begin{remark}\normalfont
	Before proceeding, we address orientability, which affects the integration framework in the volume formula. If the ambient spacetime \(M\) is not orientable, the leaves \(S_t\) may be non-orientable. In that case, integrals must be taken with respect to the canonical Riemannian density induced by \(g\) on each leaf \(S_t\), rather than a volume form (see, e.g.,~\cite[Prop.~16.45]{lee2}). The identity between the pulled-back densities remains valid without orientability or compactness, and~\eqref{Vol1} then holds as an equality of extended volumes in $[0,\infty]$. Differentiation under the integral sign and the definition of $\Theta$ below require the relevant integrals to be finite; this is ensured here by compactness.
\end{remark}

\begin{remark}\label{remark_reversing}\normalfont
	(a) A non-vanishing rescaling \(\overline U=fU\) preserves the lightlike generators and, by~\Cref{corollary:invariance_omega_en}, preserves \(\Omega\) at corresponding points. It need not preserve the constant-parameter leaves. Indeed, if
	\[
		\overline\Phi(p,s)=\Phi\bigl(p,h_p(s)\bigr),
	\]
	then
	\[
		\overline S_s=\{\Phi(p,h_p(s)):p\in S\},
	\]
	which need not coincide with any \(S_t\) when \(h_p\) depends on \(p\). Thus, the volume formula~\eqref{Vol1} is attached to the chosen geodesic extension \(U\) and its affine parameter. For a rescaled flow it retains the same form, but the parameterized family of leaves and its volume function generally change.
	
	(b) Reversing the orientation of the lightlike vector field in the radical direction (\(U \mapsto -U\)) corresponds to reversing the parametrization along each lightlike generator. As a result, whether the volume of the compact leaves is interpreted as expanding or contracting as the parameter \(t\) increases depends on the chosen time-orientation of \(U\); that is, on whether the lightlike vector field is considered future-directed or past-directed.
\end{remark}

\begin{remark} \label{remark:raychaudhuri_connection} \normalfont
	As seen in~\eqref{Vol1}, the expansion scalar $\theta$ governs the pointwise evolution of the volume density. If $\theta$ has a fixed sign on a leaf, that sign determines the sign of the corresponding volume variation; in general, the global sign is governed by the average $\Theta$ introduced below. This is consistent with the Raychaudhuri equation~\cite{beem, HawkingEllis, SanchWu}. Formula~\eqref{Vol1} provides global information for the compact leaves $S_t$ relative to the chosen future-directed geodesic generator.
\end{remark}

	Building upon these global insights, we now proceed to derive explicit formulas characterizing this evolution in the umbilical setting.
	
Since the integrand in~\eqref{Vol1} depends smoothly on $(p, t)$ for $t \in D$ and $p \in S$ (in the domain of $\Phi$), we can differentiate under the integral sign with respect to $t$, which yields the instantaneous rate of change of the volume:
\begin{equation}\label{eq:volume_derivative}
	\frac{d}{dt} \operatorname{Vol}(S_t) = \int_S \theta\bigl(\Phi(p, t)\bigr) \exp\left( \int_0^t \theta\bigl(\Phi(p, s)\bigr)\, ds \right) {\rm d}V.
\end{equation}
Both~\eqref{Vol1} and~\eqref{eq:volume_derivative} refer to the fixed geodesic extension $U$ and its affine parameter. To encode this global volume behavior, we introduce the following definition.
	
\begin{definition}\normalfont \label{def:average_expansion_scalar}
	For each $t \in D$, the \emph{average lightlike expansion scalar} $\Theta(t)$ associated with the compact leaf $S_t$ is defined by
	\[
	\Theta(t) = \frac{\int_S \theta\bigl(\Phi(p, t)\bigr)\, \exp\left( \int_0^t \theta\bigl(\Phi(p, s)\bigr)\, ds \right) \, {\rm d}V}{\int_S \exp\left( \int_0^t \theta\bigl(\Phi(p, s)\bigr)\, ds \right) \, {\rm d}V},
	\]
	where $S_t$ is a leaf of the screen distribution $\mathcal{S}$ on the lightlike hypersurface $\Sigma$ emanating from $S$ via the map $\Phi$ (as defined in~\Cref{factors}), which is assumed to be totally umbilical, and ${\rm d}V$ is the induced volume element on $S$. Here, $\theta$ denotes the lightlike expansion scalar associated with the geodesic lightlike extension $U$ of the lightlike normal vector field $\xi$ from $S$ to $\Sigma$.
	Equivalently,
	\[
		\Theta(t)=\frac{1}{\operatorname{Vol}(S_t)}
		\int_{S_t}\theta_U\,{\rm d}V_t.
	\]
	
	This scalar depends smoothly on $t \in D$ and represents a weighted average of the expansion scalar $\theta$ over the leaf $S_t$, with weight given by the accumulated volume scaling factor $\exp \bigl( \int_0^t \theta(\Phi(p, s))\, ds \bigr)$. It is defined relative to the fixed geodesic extension $U$ and its affine parameter.
\end{definition}

\begin{remark}\label{remark:Theta_conformal_dependence} \normalfont
	Let $\Theta^*$ denote the quantity obtained after a pointwise conformal change $g^*=e^{2u}g$, using the $g^*$-geodesic extension with the same initial lightlike normal field. In general, neither the resulting parametrized family of leaves nor $\Theta$ is preserved. Both the induced volume elements and the conformal factor relating corresponding points are affected; the latter obeys the covariance law in~\Cref{remark:Omega_conformal_invariance}. If $u$ is constant, however, the Levi-Civita connection and the affine parametrization are unchanged, while all induced volume elements acquire the same constant multiple. In that case the quotient defining $\Theta(t)$ is unchanged.
\end{remark}

	With this definition, we can now state the following result that determines the volume evolution of compact spacelike leaves $S_t$ of the foliation:
	
	\begin{theorem}\label{propo:volume evolution formula}
		Under the assumptions of~\Cref{Th:Vol1}, the volume of each compact leaf $S_t$ of the screen distribution ${\mathcal S}$ evolves according to the formula:
		\begin{equation}\label{volume evolution formula}
			\operatorname{Vol}(S_t) = \operatorname{Vol}(S) \exp\left( \int_0^t \Theta(s) \,ds \right).
		\end{equation}
		This formula holds for all $t \in D$ for the fixed geodesic extension $U$ and its affine parameter.
	\end{theorem}
\begin{proof}
Taking into account~\eqref{Vol1} and \eqref{eq:volume_derivative}, the definition of \(\Theta(t)\) leads to the ordinary differential equation
	\[
	\frac{d}{dt} \operatorname{Vol}(S_t) = \Theta(t) \operatorname{Vol}(S_t),
	\]
	whose solution, with initial condition $\operatorname{Vol}(S_0) = \operatorname{Vol}(S)$, gives the stated formula~\eqref{volume evolution formula}.
\end{proof}

\begin{remark}\normalfont
	The average lightlike expansion scalar \(\Theta(t)\) and the resulting volume evolution formula~\eqref{volume evolution formula} provide a powerful and geometrically intuitive framework for investigating the global volume evolution of the compact spacelike leaves \(S_t\) along the lightlike umbilical direction (for \(t \in D\)). The sign of the average lightlike expansion scalar \(\Theta(t)\) precisely dictates the instantaneous rate of change of the leaf volume \(S_t\): a positive (respectively negative) value indicates that the volume is instantaneously increasing (respectively decreasing), reflecting a dominance of expansion (respectively contraction) of the corresponding lightlike generators. Moreover, the integral \(\int_0^t \Theta(s) \, ds\) in~\eqref{volume evolution formula} elegantly governs the \emph{overall} volume scaling of \(S_t\) relative to the initial submanifold \(S\), capturing the accumulated effect of expansion and contraction along the lightlike flow.
	
	Crucially, the total umbilicity of the generated lightlike hypersurface \(\Sigma\) leads to an underlying geometric simplicity that allows for a direct and transparent connection between the average expansion scalar \(\Theta(t)\) and the global volume evolution of the leaves.
	
	This parallels the spirit of Hawking's Black Hole Area Theorem~\cite{HawkingEllis, Wald}: relative to the chosen future-directed geodesic generator, if the average lightlike expansion is non-negative, the volume of the leaves does not decrease. Reversing the time-orientation reverses the corresponding interpretation of increase and decrease, as discussed in~\Cref{remark_reversing}(b).
\end{remark}

	\subsection{Variational Properties}
	
	Globally defined flow-invariant vector fields on~$\Sigma$, when restricted to a lightlike geodesic generator~$\gamma$, yield $S$-Jacobi fields along~$\gamma$. We refer the reader to~\cite{beem, oneill} for standard definitions and properties of Jacobi fields along lightlike geodesics.
	
	\begin{theorem}\label{corSJacobi}
		Let $S$ be a codimension-two spacelike submanifold embedded in a spacetime $M$, and let $\xi$ be a lightlike normal vector field on $S$. Let $\Sigma$ be the lightlike hypersurface generated from $S$ by the normal exponential map along $\xi$, as in~\Cref{factors}. 
		
		Then, for any flow-invariant vector field $X \in \mathfrak{I}(\Sigma)$ and any lightlike geodesic generator $\gamma$ of $\Sigma$, the vector field $J = X \circ \gamma$ is an $S$-Jacobi field along $\gamma$, which describes the variation of $\gamma$ through lightlike geodesics initially normal to $S$.
		
		Moreover, if $\Sigma$ is totally umbilical, then $J= X \circ \gamma$ satisfies the second-order differential equation
		\begin{equation}\label{JRicc}
			J'' + \frac{\widetilde{{\rm Ric}}(\gamma', \gamma')}{n} J = f\, \gamma',
		\end{equation}
		where $n = \dim S$, and $\widetilde{{\rm Ric}}$ is the Ricci tensor of $M$. Here, $f$ is a smooth function along $\gamma$ given explicitly by
		\begin{equation}\label{eq:f_expresion_Variational_Prp}
			f = \gamma'(\tau(J)) - \mu\, \tau(J),
		\end{equation}
		where $\mu$ is the umbilicity function and $\tau$ is the rotation one-form associated to the geodesic extension of $\xi$ to $\Sigma$ and to the screen distribution $\mathcal{S}$, both canonically determined by~\Cref{factors}.
		
		Primes denote covariant derivatives along $\gamma$, that is, $J' = \widetilde{\nabla}_{\gamma'} J$ and $J'' = \widetilde{\nabla}_{\gamma'} \widetilde{\nabla}_{\gamma'} J$.
	\end{theorem}
	\begin{proof}
		Let $U$ be the geodesic extension of $\xi$ to a lightlike vector field on~$\Sigma$, and let $\gamma$ be a lightlike geodesic generator of $\Sigma$ with $\gamma(0) \in S$ and $\gamma' = U$ along~$\gamma$.
		
		Let $X \in \mathfrak{I}(\Sigma)$ be a spacelike vector field, where $\mathfrak{I}(\Sigma)$ denotes the set of flow-invariant vector fields on $\Sigma$ obtained as pushforwards of horizontal lifts of vector fields on $S$. The flow-invariance of $X$ implies
		\begin{equation}\label{eq:flow_invariance}
			\widetilde{\nabla}_U X = \widetilde{\nabla}_X U.
		\end{equation}
		Since $U$ is geodesic, that is, $\widetilde{\nabla}_U U = 0$, the curvature tensor of $M$ satisfies
		\begin{align}
			\widetilde{R}(X, U)U &= \widetilde{\nabla}_X \widetilde{\nabla}_U U - \widetilde{\nabla}_U \widetilde{\nabla}_X U - \widetilde{\nabla}_{[X, U]} U \nonumber \\
			&= -\widetilde{\nabla}_U \widetilde{\nabla}_X U \label{Eq:R_JacobiOp} \\
			&= -\widetilde{\nabla}_U \widetilde{\nabla}_U X. \nonumber
		\end{align}
		This confirms that $J'' + \widetilde{R}(J, \gamma')\gamma' = 0$ for $J = X \circ \gamma$, establishing that $X \circ \gamma$ is a Jacobi field along~$\gamma$.
		
		The initial conditions for an $S$-Jacobi field, namely
		\begin{center}
			$J(0) \in T_{\gamma(0)}S$\quad and \quad$\widetilde{\nabla}_{\gamma'(0)} J + A_{\gamma'(0)}J(0) \perp T_{\gamma(0)}S$,
		\end{center}
		are satisfied by construction of $J$, the Weingarten formula~\eqref{fWeingarten}, and the commuting property~\eqref{eq:flow_invariance}. Consequently, $J = X \circ \gamma$ is an $S$-Jacobi field, representing the variation vector field of a variation of~$\gamma$ through lightlike geodesics normal to $S$, as $g(J,\gamma')=0$~\cite[Cor.~10.40]{oneill}.
		
		Now, suppose $\Sigma$ is totally umbilical, characterized by a lightlike second fundamental form $B_U = \mu g$ for some smooth function $\mu$ on $\Sigma$. Making use of~\eqref{RelaLightlik}, the equation~\eqref{LightlikeDes1} then takes the form
		\begin{equation}
			\widetilde{\nabla}_X U = -\mu X + \tau(X) U, \label{Eq:nablaXU}
		\end{equation}
		for any $X \in \Gamma({\mathcal S})$, where $\tau$ is the rotation one-form associated with $U$ and the screen distribution $\mathcal{S}$.
		
		Next, consider the Jacobi operator $\widetilde{R}_U \colon \mathfrak{X}(\Sigma) \to \mathfrak{X}(\Sigma)$ (see~\cite[p.~219]{oneill}), defined by
		\[
		\widetilde{R}_U(X) = \widetilde{R}(X, U)U.
		\]
		Applying~\eqref{Eq:R_JacobiOp} and~\eqref{Eq:nablaXU}, we find that for any $X, Y \in \mathfrak{I}(\Sigma)$,
		\begin{align}
			g(\widetilde{R}_U(X), Y) 
			&= g(-\widetilde{\nabla}_U \widetilde{\nabla}_X U, Y) \nonumber \\
			&= g(-\widetilde{\nabla}_U (-\mu X + \tau(X) U), Y) \nonumber \\
			&= (U(\mu) - \mu^2) g(X, Y), \label{Eq:Inv_Cur}
		\end{align}
		where we have used $g(U, Y) = 0$ and $\widetilde{\nabla}_U U = 0$. Since the vector fields in $\mathfrak{I}(\Sigma)$ span the screen distribution~$\mathcal{S}$ at each point of $\Sigma$, equation~\eqref{Eq:Inv_Cur} holds for all $X, Y \in \Gamma(\mathcal{S})$. 
		
		To extend this identity to arbitrary vector fields in $\mathfrak{X}(\Sigma)$, note that any such vector field can be written as the sum of a section of $\mathcal{S}$ and a vector field collinear with $U$. By linearity and the properties of the curvature tensor, it follows from~\eqref{Eq:Inv_Cur} that
		\begin{equation}\label{Eq:Jacobi_operator_Umbilical}
			g(\widetilde{R}_U(X), Y) = (U(\mu) - \mu^2) g(X, Y),
		\end{equation}
		for all $X, Y \in \mathfrak{X}(\Sigma)$.
		
		The Ricci tensor component $\widetilde{{\rm Ric}}(U, U)$ is related to the Jacobi operator trace through $\widetilde{{\rm Ric}}(U, U) =  \mathrm{trace}(\widetilde{R}_U)$ (see~\cite[Lem.~8.9]{oneill}). Since $n = \dim S$, it follows from~\eqref{Eq:Jacobi_operator_Umbilical} that $\widetilde{{\rm Ric}}(U, U) = n(U(\mu) - \mu^2)$. Therefore, the Jacobi equation for a spacelike $S$-Jacobi field of the form $J = X \circ \gamma$ with $X \in \mathfrak{I}(\Sigma)$ can be written as~\eqref{JRicc}.
		
		Finally, an explicit expression for $f$ can be computed as follows. Along the geodesic generator $\gamma$, the flow-invariance property~\eqref{eq:flow_invariance} and the Weingarten equation~\eqref{Eq:nablaXU} give $J' = -\mu J + \tau(J) \gamma'$. Differentiating and using $\widetilde{\nabla}_{\gamma'} \gamma' = 0$ yields $J'' = -\gamma'(\mu) J - \mu J' + \gamma'(\tau(J))\, \gamma'$. Substituting the expression for $J'$ into this equation, we find
		\begin{equation*}
			J'' + (\gamma'(\mu) - \mu^2) J = (\gamma'(\tau(J)) - \mu \tau(J))\, \gamma'.
		\end{equation*}
		Therefore, the explicit formula for $f$ is as given in~\eqref{eq:f_expresion_Variational_Prp}.
	\end{proof}
	
	\begin{corollary}\label{cor:umbilical-jacobi}
		In the totally umbilical case of~\Cref{corSJacobi}, every Jacobi field $J$ along a lightlike geodesic generator $\gamma$ of $\Sigma$ with $g(J, \gamma') = 0$ satisfies
		\begin{equation}\label{JRicc_2}
			J'' + \frac{\widetilde{\mathrm{Ric}}(\gamma',\gamma')}{n} J = f\, \gamma',
		\end{equation}
		for some smooth function $f$ along $\gamma$. In this setting, $J$ does not need to be tangent to the screen distribution $\mathcal{S}$, and therefore $f$ is not, in general, related to the associated rotation one-form $\tau$, nor can it necessarily be written in the explicit form given in~\eqref{eq:f_expresion_Variational_Prp}.
	\end{corollary}
	\begin{proof}
		Since equation~\eqref{Eq:Jacobi_operator_Umbilical} holds for all vector fields $X, Y$ tangent to $\Sigma$, it follows that in the umbilical case, every Jacobi field $J$ along a lightlike geodesic generator $\gamma$ of~$\Sigma$ with $g(J, \gamma') = 0$ satisfies the stated equation \eqref{JRicc_2}.
	\end{proof}

\begin{remark}\normalfont
	When $\Sigma$ is totally umbilical, the Jacobi equation~\eqref{JRicc_2} takes on a distinctive form when considered in the appropriate quotient space of vector fields along the lightlike geodesic generators of~$\Sigma$ (see~\cite{beem, MinguzziSanchez} for details). This perspective naturally motivates the introduction of the concept of \emph{Jacobi classes} for a more refined analysis, especially in the study of \emph{quotient focal points} along those lightlike geodesics~\cite{beem}, which, by construction, emanate orthogonally from the initial codimension-two spacelike submanifold~$S$.
\end{remark}
	
		\begin{remark}\normalfont
		Theorem~\ref{corSJacobi} extends \cite[Prop.~3.6]{gutierrez15}, which analyzes Jacobi fields in totally umbilical light cones of generalized Robertson--Walker spacetimes, where the rotation one-form $\tau$ vanishes as a consequence of the geometry. The case $\tau = 0$ is especially interesting, as the Jacobi equation~\eqref{JRicc} simplifies considerably through~\eqref{eq:f_expresion_Variational_Prp}. In contrast, our result generalizes this setting by considering lightlike geodesics emanating from a general codimension-two spacelike submanifold $S$ in more general spacetimes, thereby allowing for a non-zero $\tau$ that directly influences the umbilical Jacobi equation~\eqref{JRicc}. This generalization is particularly relevant in situations where the normal connection of each leaf is not flat, as will be discussed in the following subsection on parallelism.
	\end{remark}

\subsection{Curvature and Parallelism}

We now apply the characterization of parallel umbilical lightlike normal sections from~\Cref{sec:para} to totally umbilical lightlike hypersurfaces and their associated integrable screen distributions, leading to the following result.

\begin{theorem}\label{corollary:parallel_normal_connection}
	Let $S$ be a codimension-two spacelike submanifold embedded in a spacetime $M$, and suppose that $S$ admits a lightlike normal vector field $\xi$. Let $\Sigma$ denote the lightlike hypersurface emanating from $S$, constructed as in~\Cref{factors}, and assume explicitly that $\Sigma$ is totally umbilical. Let $\mathcal{S}$ be the associated integrable screen distribution on $\Sigma$, and let $U$ be a lightlike vector field on $\Sigma$ with $V$ the transverse lightlike vector field normalized so that $g(U,V) = -1$. Then, the following properties hold:
	\begin{enumerate}
		\item[(a)] The curvature tensor $\widetilde{R}$ of $M$ satisfies
		\begin{equation}\label{eq:curvature_tau}
			g(\widetilde{R}(X,Y)U,V) = -d\tau(X,Y),
		\end{equation}
		for all $X, Y \in \Gamma(\mathcal{S})$, where $\tau$ is the rotation one-form associated with $U$ and $\mathcal{S}$. Moreover, the expression in~\eqref{eq:curvature_tau} is invariant under conformal changes of the Lorentzian metric of $M$.
		
		\item[(b)] The vector field $U$ can be locally rescaled along each leaf of the screen distribution $\mathcal{S}$ so that it is parallel (with respect to the induced normal connection) if and only if the curvature tensor of $M$ satisfies
		\begin{equation}\label{eq:vanishing_curvature}
			g(\widetilde{R}(X,Y)U,V) = 0,
		\end{equation}
		for all $X, Y \in \Gamma(\mathcal{S})$. Under this condition, the normal connection of each leaf of $\mathcal{S}$ is flat. If, in addition, $S$ is simply connected, then the normal holonomy group of every globally defined leaf $S_t$, $t\in D$, is trivial. Moreover, this property, namely the local rescalability of $U$ to a parallel vector field, is also conformally invariant.
		
		\item[(c)] Condition~\eqref{eq:vanishing_curvature} holds whenever the spacetime $M$ is locally conformally flat.
	\end{enumerate}
\end{theorem}

\begin{proof}
	\emph{(a)} For each leaf of the integrable screen distribution $\mathcal{S}$, the restriction of $U$ defines a lightlike normal vector field, which is umbilical by~\Cref{cor:umbilical_leaves}. The compatibility between the Weingarten formula for codimension-two spacelike submanifolds~\eqref{fWeingarten} and that for the lightlike hypersurface~\eqref{LightlikeDes1} ensures that the rotation one-form $\tau$ restricts naturally to each leaf. The curvature formula~\eqref{eq:RNorTau} together with the Ricci equation~\eqref{eq:RSpacetimePerp} then yields $g(\widetilde{R}(X,Y)U,V) = -d\tau(X,Y)$ for all $X, Y \in \Gamma(\mathcal{S})$, since $\mathcal{S}$ is integrable. 
	
The conformal invariance of this expression follows from the fact that the normal curvature tensor associated to each leaf of the screen distribution is conformally invariant under changes of the ambient Lorentzian metric. Since $d\tau$ is related to the normal curvature tensor by equation~\eqref{eq:RNorTau}, it follows that $d\tau$ is also conformally invariant. Therefore, equation~\eqref{eq:curvature_tau} is preserved under conformal transformations.
	
	\emph{(b)} The equivalence follows by applying~\Cref{CriterionParaUmbi} to each leaf of $\mathcal{S}$: for any leaf, the restriction of $U$ defines an umbilical lightlike normal vector field, and $U$ can be locally rescaled along the leaf to a parallel vector field with respect to the induced normal connection if and only if $g(\widetilde{R}(X,Y)U,V) = 0$ for all $X, Y \in \Gamma(\mathcal{S})$. In this case, the normal connection of each leaf is flat. If $S$ is simply connected, then every globally defined leaf $S_t$, $t\in D$, is simply connected because $\Phi_t\colon S\to S_t$ is a diffeomorphism; its normal holonomy group is therefore trivial by~\Cref{PUS}. The conformal invariance of the local rescalability of $U$ to a parallel vector field along each leaf follows from the fact that the condition $d\tau = 0$ is conformally invariant, as explained above.
	
	\emph{(c)} In a locally conformally flat spacetime, each leaf of $\mathcal{S}$ has flat normal connection, and condition~\eqref{eq:vanishing_curvature} holds.
\end{proof}

\begin{remark}\normalfont
	The curvature property~\eqref{eq:curvature_tau} is consistent with the independence of the exterior derivative $d\tau$ under rescalings of the lightlike vector field $U \in \mathfrak{X}(\Sigma)$, as discussed in~\Cref{remark:dtau_ditinghished}. It should be emphasized that the lightlike vector field $U$ in~\Cref{corollary:parallel_normal_connection} is not required to be the geodesic extension of $\xi$ to $\Sigma$ constructed in~\Cref{factors}; in fact, $U$ can be any lightlike vector field generating the radical distribution of $\Sigma$.
\end{remark}

\begin{remark}\normalfont
Condition~\eqref{eq:vanishing_curvature} holds on the event horizon of
the Kruskal spacetime; see~\cite[p.~100]{kupeli}.
\end{remark}

	\section*{Conclusion}

	In summary, every embedded codimension-two spacelike submanifold in a spacetime carrying a lightlike normal section admits the lightlike factorization described in~\Cref{factors}. The associated lightlike hypersurface is totally umbilical precisely when the restriction of its lightlike generator is an umbilical normal section along every leaf of the integrable screen distribution. In particular, the umbilicity of the initial lightlike normal section suffices when the ambient spacetime is locally conformally flat. Under the explicit total umbilicity hypothesis, the conformal, volume, variational, and curvature results of~\Cref{sec:conse_Umbil_Fac_Th} apply under the hypotheses stated there.
	
	\section*{Acknowledgements}
	I would like to thank Francisco J. Palomo (University of Málaga) and Alfonso Romero (University of Granada) for their detailed reviews and valuable feedback, which have significantly improved this paper. 
	
	I also thank the anonymous referees for their careful reading and constructive comments, which have helped to further improve the clarity and quality of the manuscript.
	

\end{document}